\documentclass[12pt,twoside,leqno]{article}
\usepackage{amsmath}
\usepackage{amssymb}
\usepackage{amsxtra}
\usepackage{amscd}
\usepackage{amsthm}
\usepackage[mathscr]{eucal}
 \usepackage{oldgerm,extarrows,color,graphicx,hyperref} 

\usepackage[nohug]{diagrams}

\theoremstyle{plain}
\newtheorem{thm}[subsection]{Theorem}
\newtheorem{prop}[subsection]{Proposition}
\newtheorem{cor}[subsection]{Corollary}
\newtheorem{lem}[subsection]{Lemma}

\theoremstyle{definition}

\newtheorem{rem}[subsection]{Remark}
\newtheorem{para}[subsection]{}

\newtheorem{sbpara}[subsubsection]{}

\newenvironment{pf}{\proof[\proofname]}{\endproof}
\newenvironment{pf*}[1]{\proof[#1]}{\endproof}

\begin{document}

\title{Syntomic complexes of $F$-crystals and Tamagawa number conjecture  in characteristic $p$, by Kazuya Kato,
with Appendix by Arthur Ogus}


\maketitle

\newcommand\Cal{\mathcal}
\newcommand\define{\newcommand}

\define\gp{\mathrm{gp}}%
\define\fs{\mathrm{fs}}%
\define\an{\mathrm{an}}%
\define\mult{\mathrm{mult}}%
\define\Ker{\mathrm{Ker}\,}%
\define\Coker{\mathrm{Coker}\,}%
\define\Hom{\mathrm{Hom}\,}%
\define\Ext{\mathrm{Ext}\,}%
\define\rank{\mathrm{rank}\,}%
\define\gr{\mathrm{gr}}%
\define\cHom{\Cal{Hom}}
\define\cExt{\Cal Ext\,}%

\define\cC{\Cal C}
\define\cD{\Cal D}
\define\cO{\Cal O}
\define\cS{\Cal S}
\define\cM{\Cal M}
\define\cG{\Cal G}
\define\cH{\Cal H}
\define\cE{\Cal E}
\define\cF{\Cal F}
\define\cN{\Cal N}
\define\cT{\Cal T}
\define\cP{\Cal P}
\define\cQ{\Cal Q}
\define\fF{\frak F}
\define\Dc{\check{D}}
\define\Ec{\check{E}}

\newcommand{\N}{{\mathbb{N}}}
\newcommand{\Q}{{\mathbb{Q}}}
\newcommand{\Z}{{\mathbb{Z}}}
\newcommand{\R}{{\mathbb{R}}}
\newcommand{\C}{{\mathbb{C}}}
\newcommand{\bN}{{\mathbb{N}}}
\newcommand{\bQ}{{\mathbb{Q}}}
\newcommand{\bF}{{\mathbb{F}}}
\newcommand{\bZ}{{\mathbb{Z}}}
\newcommand{\bP}{{\mathbb{P}}}
\newcommand{\bR}{{\mathbb{R}}}
\newcommand{\bC}{{\mathbb{C}}}
\newcommand{\bbQ}{{\bar \mathbb{Q}}}
\newcommand{\ol}[1]{\overline{#1}}
\newcommand{\too}{\longrightarrow}
\newcommand{\respect}{\rightsquigarrow}
\newcommand{\compatible}{\leftrightsquigarrow}
\newcommand{\upc}[1]{\overset {\lower 0.3ex \hbox{${\;}_{\circ}$}}{#1}}
\newcommand{\Gmlog}{\bG_{m, \log}}
\newcommand{\Gm}{\bG_m}
\newcommand{\ep}{\varepsilon}
\newcommand{\Spec}{\operatorname{Spec}}
\newcommand{\val}{{\mathrm{val}}} 
\newcommand{\n}{\operatorname{naive}}
\newcommand{\bs}{\operatorname{\backslash}}
\newcommand{\Gal}{\operatorname{{Gal}}}
\newcommand{\gal}{{\rm {Gal}}({\bar \Q}/{\Q})}
\newcommand{\galp}{{\rm {Gal}}({\bar \Q}_p/{\Q}_p)}
\newcommand{\gall}{{\rm{Gal}}({\bar \Q}_\ell/\Q_\ell)}
\newcommand{\wep}{W({\bar \Q}_p/\Q_p)}
\newcommand{\wel}{W({\bar \Q}_\ell/\Q_\ell)}
\newcommand{\Ad}{{\rm{Ad}}}
\newcommand{\BS}{{\rm {BS}}}
\newcommand{\even}{\operatorname{even}}
\newcommand{\End}{{\rm {End}}}
\newcommand{\odd}{\operatorname{odd}}
\newcommand{\GL}{\operatorname{GL}}
\newcommand{\np}{\text{non-$p$}}
\newcommand{\g}{{\gamma}}
\newcommand{\G}{{\Gamma}}
\newcommand{\Lam}{{\Lambda}}
\newcommand{\La}{{\Lambda}}
\newcommand{\lam}{{\lambda}}
\newcommand{\la}{{\lambda}}
\newcommand{\uL}{{{\hat {L}}^{\rm {ur}}}}
\newcommand{\uQp}{{{\hat \Q}_p}^{\text{ur}}}
\newcommand{\sel}{\operatorname{Sel}}
\newcommand{\dt}{{\rm{Det}}}
\newcommand{\Sig}{\Sigma}
\newcommand{\fil}{{\rm{fil}}}
\newcommand{\SL}{{\rm{SL}}}
\newcommand{\spl}{{\rm{spl}}}
\newcommand{\st}{{\rm{st}}}
\newcommand{\Isom}{{\rm {Isom}}}
\newcommand{\Mor}{{\rm {Mor}}}
\newcommand{\bg}{\bar{g}}
\newcommand{\id}{{\rm {id}}}
\newcommand{\cone}{{\rm {cone}}}
\newcommand{\al}{a}
\newcommand{\ChL}{{\cal{C}}(\La)}
\newcommand{\Image}{{\rm {Image}}}
\newcommand{\toric}{{\operatorname{toric}}}
\newcommand{\torus}{{\operatorname{torus}}}
\newcommand{\Aut}{{\rm {Aut}}}
\newcommand{\Qp}{{\mathbb{Q}}_p}
\newcommand{\barQp}{{\mathbb{Q}}_p}
\newcommand{\Qpur}{{\mathbb{Q}}_p^{\rm {ur}}}
\newcommand{\Zp}{{\mathbb{Z}}_p}
\newcommand{\Zl}{{\mathbb{Z}}_l}
\newcommand{\Ql}{{\mathbb{Q}}_l}
\newcommand{\Qlur}{{\mathbb{Q}}_l^{\rm {ur}}}
\newcommand{\F}{{\mathbb{F}}}
\newcommand{\eps}{{\epsilon}}
\newcommand{\epsLa}{{\epsilon}_{\La}}
\newcommand{\epsLaVxi}{{\epsilon}_{\La}(V, \xi)}
\newcommand{\epsOLaVxi}{{\epsilon}_{0,\La}(V, \xi)}
\newcommand{\Qplin}{{\mathbb{Q}}_p(\mu_{l^{\infty}})}
\newcommand{\otimesQplin}{\otimes_{\Qp}{\mathbb{Q}}_p(\mu_{l^{\infty}})}
\newcommand{\galFl}{{\rm{Gal}}({\bar {\Bbb F}}_\ell/{\Bbb F}_\ell)}
\newcommand{\gallur}{{\rm{Gal}}({\bar \Q}_\ell/\Q_\ell^{\rm {ur}})}
\newcommand{\galFF}{{\rm {Gal}}(F_{\infty}/F)}
\newcommand{\galFv}{{\rm {Gal}}(\bar{F}_v/F_v)}
\newcommand{\galF}{{\rm {Gal}}(\bar{F}/F)}
\newcommand{\epsVxi}{{\epsilon}(V, \xi)}
\newcommand{\epsOVxi}{{\epsilon}_0(V, \xi)}
\newcommand{\plim}{\lim_
{\scriptstyle 
\longleftarrow \atop \scriptstyle n}}
\newcommand{\sig}{{\sigma}}
\newcommand{\ga}{{\gamma}}
\newcommand{\del}{{\delta}}
\newcommand{\Vss}{V^{\rm {ss}}}
\newcommand{\Bst}{B_{\rm {st}}}
\newcommand{\Dpst}{D_{\rm {pst}}}
\newcommand{\Dcrys}{D_{\rm {crys}}}
\newcommand{\DdR}{D_{\rm {dR}}}
\newcommand{\Fin}{F_{\infty}}
\newcommand{\Kla}{K_{\lambda}}
\newcommand{\Ola}{O_{\lambda}}
\newcommand{\Mla}{M_{\lambda}}
\newcommand{\Det}{{\rm{Det}}}
\newcommand{\Sym}{{\rm{Sym}}}
\newcommand{\LaSa}{{\La_{S^*}}}
\newcommand{\cX}{{\cal {X}}}
\newcommand{\MHG}{{\frak {M}}_H(G)}
\newcommand{\tauMla}{\tau(M_{\lambda})}
\newcommand{\Fvur}{{F_v^{\rm {ur}}}}
\newcommand{\Lie}{{\rm {Lie}}}
\newcommand{\cL}{{\cal {L}}}
\newcommand{\cW}{{\cal {W}}}
\newcommand{\fq}{{\frak {q}}}
\newcommand{\cont}{{\rm {cont}}}
\newcommand{\SC}{{SC}}
\newcommand{\Om}{{\Omega}}
\newcommand{\dR}{{\rm {dR}}}
\newcommand{\crys}{{\rm {crys}}}
\newcommand{\hatSig}{{\hat{\Sigma}}}
\newcommand{\rdet}{{{\rm {det}}}}
\newcommand{\ord}{{{\rm {ord}}}}
\newcommand{\BdR}{{B_{\rm {dR}}}}
\newcommand{\BdRO}{{B^0_{\rm {dR}}}}
\newcommand{\Bcrys}{{B_{\rm {crys}}}}
\newcommand{\Qw}{{\mathbb{Q}}_w}
\newcommand{\barkappa}{{\bar{\kappa}}}
\newcommand{\cZ}{{\Cal {Z}}}
\newcommand{\oppLa}{{\Lambda^{\circ}}}
\newcommand{\bG}{{\mathbb{G}}}
\newcommand{\Rep}{{{\rm Rep}}}
\newcommand{\red}{{{\rm red}}}
\newcommand{\DR}{{{\rm DR}}}
\newcommand{\tE}{{\tilde E}}
\newcommand{\tnabla}{{\tilde \nabla}}
\newcommand{\tPhi}{{\tilde \Phi}}
\newcommand{\cx}{\cdot}
\newcommand{\ot}{\otimes}

\begin{abstract} 
We consider syntomic complexes of uniform $F$-crystals and relate them 
to Tamagawa number conjecture in characteristic $p$.

\end{abstract}

{\bf Contents}

\medskip

\S\ref{s:0}. Introduction.

\S\ref{s:Fcr}. Uniform $F$-crystals.

\S\ref{s:DR}. De Rham complexes.

\S\ref{s:Syn}. Syntomic complexes. 

\S\ref{s:Tam}. On Tamagawa number conjecture.

\S\ref{s:Cpl}. Complements. 

 \S\ref{s:App}. Appendix: Remarks on F-spans

\section{Introduction}\label{s:0}
Uniform $F$-crystals in the sense of Ogus have good integral properties, good notion of exact sequences,  and form a nice category which is stable under tensor products and  duals.  In this paper, we define syntomic complexes of uniform $F$-crystals. We prove that the group of extension classes of uniform $F$-crystals is expressed by the first cohomology of a syntomic complex (Theorem \ref{H1sD}). By using the syntomic complexes, we consider the Tamagawa number conjecture in characteristic $p$, which is an analogue of Tamagawa number conjecture ([BK]) of S. Bloch and the author over number fields.

In Section 2 and Section 3, we review the parts on uniform $F$-crystals and their modified de Rham complexes  in the works \cite{O0} and \cite{O} of Ogus. We present the theory in the style which is convenient for the descriptions of Sections 4--6.  In Section 4, we define the syntomic complex of a uniform $F$-crystal by using these modified de Rham complexes. In \S5, we consider Tamagawa number conjecture in characteristic $p$. In Section 6, we give comments, especially, on a relation with the theory of Frobenius gauges of J.-M. Fontaine and U. Janssen.

The syntomic complexes of F-crystals in this paper are considered 
 independently in the paper \cite{BMS} of B. Bhatt, M. Morrow and P. Scholze and in 
the work of Bhatt and Scholze on prismatic cohomology. Tamagawa number conjecture in characteristic $p$ was also considered recently by O. Brinon and F. Trihan \cite{BT} by a different method. 

In Appendix, A. Ogus improves Proposition \ref{filD} of this paper and explains the relation to the work \cite{Sh} of A. Shiho.  (See  Remark \ref{Onew}). 

The author expresses his hearty thanks to  Arthur Ogus for valuable advice. He is also thankful to Teruhisa Koshikawa for the help.

\section{Uniform $F$-crystals}\label{s:Fcr}

In Section 2 and Section 3, we review the parts of the works of Ogus (\cite{O0}, \cite{O}) concerning  uniform $F$-crystals. 

 Uniform $F$-crystals have good integral properties. The category of uniform $F$-crystals is stable under direct sum, tensor products, duals, and Tate twists. Nice notion of exact sequence of uniform $F$-crystals is defined.  
 
 The book \cite{O} is written in big generality, and the author thinks that it is nice to pick up this important part on uniform $F$-crystals, and give a short self-contained exposition.  The author sometimes present proofs  which are different from the original and which seem to fit uniform $F$-crystals well. 
In Section 2 and Section 3, the only new results which are not contained in \cite{O0} and \cite{O} are  \ref{filD}, \ref{Hodge}, \ref{DRrr2}, in which we remove the assumption in \cite{O0} and \cite{O} that the lifting of the Frobenius is parallelizable in the sense of \cite{O} 1.2.6  (but see Remark \ref{Onew}).

\begin{para}\label{XD} Let $k$ be a perfect field  of characteristic $p>0$ and 
let $X$ be a log smooth fs log scheme over $k$ (\cite{Il}).

Let $D$ be a
crystal on the log crystalline site of $X$ over $W(k)$ (\cite{KFI} Section 5) which is  locally free of finite rank,  and let $$\Phi: F^*D_{Q_p}\overset{\cong}\to D_{\Q_p}$$ be an isomorphism given in the category in which $\Hom$ is replaced by $\Q_p\otimes  \Hom$. Here $F$ denotes the Frobenius ($p$-th power morphism) $F: X\to X$ and $F^*D$ denotes the pullback of $D$ by $F$. 

We call this object $D=(D, \Phi)$ an $F$-crystal on $X$. We call $\Phi$ the Frobenius operator of $D=(D, \Phi)$. 

For $q\in \Z$, let $$\omega^q_X:=\Omega^q_X(\log)$$ be the sheaf of differential $q$-forms on $X$ with log poles (\cite{KFI}).

We denote by  $D_X$  the vector bundle on $X$ with an integrable connection
$$\nabla: D_X \to D_X\otimes_{\cO_X} \omega^1_X$$
associated to $D$. 
\end{para}

\begin{para}\label{dual} Duals. Tate twists. 

For an $F$-crystal $D=(D, \Phi)$ on $X$, its dual $D^*=(D^*, \Phi^*)$ is defined as follows. The crystal $D^*$ is the linear dual of $D$. The Frobenius operator $\Phi^*$ of $D^*$ is the inverse $({}^t\Phi)^{-1}$ of the transpose ${}^t\Phi: D_{\Q_p} \to F^*D_{\Q_p}$ of $\Phi$. 

For an $F$-crystal $D=(D, \Phi)$ on $X$ and for $r\in \Z$, its $r$-th Tate twist $D(r)$ is defined to be the pair $(D, p^{-r}\Phi)$. 

Let the $F$-crystal ${\bf 1}$ be the structure sheaf of the log crystalline site endowed with $\Phi=\text{id}$. We have
$$D(r)= D\otimes {\bf 1}(r).$$

Even if the Frobenius operator $\Phi$ of $D$  comes from  a homomorphism $F^*D\to D$ (without denominator), the Frobenius operators of the dual of $D$ and of the Tate twists of $D$
 may have denominators.  To have the nice notions of dual and Tate twists in this paper, we allow that the Frobenuis operator $\Phi$ of an $F$-crystal can have a denominator as in \ref{XD}.

\end{para}

\begin{para}\label{XZ}  By a $p$-adic formal fs log scheme over $W(k)$, we mean a family $Y=(Y_n)_{n\geq 1}$ of fs log schemes $Y_n$ over $W_n(k)$ endowed with an isomorphism $Y_{n+1}\otimes_{W_{n+1}(k)} W_n(k)\cong Y_n$ for each $n\geq 1$. (We regard $Y$ naturally  as a formal scheme in the sense of Grothendieck. So the underlying topological space of $Y$ is that of $Y_n$ for any $n\geq 1$,)
We say $Y$ is log smooth if $Y_n$ is log smooth over $W_n(k)$ for any $n\geq 1$. We regard $Y$ as a formal scheme in the sense of Grothendieck 

Assume that we have a closed immersion $X\subset Y$ over $W(k)$ where $Y$ is a log smooth $p$-adic formal  fs log scheme over $W(k)$ 
endowed with a morphism $F_Y:Y\to Y$ which lifts the Frobenius of $Y\otimes \Z/p\Z$.

We call this situation $X\subset Y$ an embedded situation. 

We have an embedded situation \'etale locally on $X$. 

Let $D_X(Y)$ be the PD envelop of $X$ in $Y$ (\cite{KFI} Section 5). It is a $p$-adic formal fs log scheme over $W(k)$, its underlying topological space is identified with that of $X$, and the log structure of $X$ coincides with the pullback of the log structure of $D_X(Y)$. Let 
$\cO_{D_X(Y)}$ be  the structure sheaf of the formal scheme $D_X(Y)$. 
Then $\cO_{D_X(Y)}$ is $p$-torsion free. Let $$\omega^q_Y:=\Omega^q_Y(\log)$$ be the sheaf of differential $q$-forms on $Y$ with log poles.

  Then the crystal $D$ corresponds  to a locally free 
$\cO_{D_X(Y)}$-module $D_Y$ of finite rank on $D_X(Y)$ endowed with an  integrable connection $$\nabla: D_Y\to D_Y\otimes_{\cO_Y} \omega^1_Y$$ whose mod $p^n$  satisfies the logarithmic version (the condition (iii) of Theorem 6.2 in \cite{KFI}) of the quasi-nilpotence in \cite{BO} 4.10, and the $F$-crystal corresponds to such $D_Y$ endowed with an isomorphism of $\Q_p\otimes_{\Z_p} \cO_{D_X(Y)}$-modules $$\Phi:\Q_p\otimes_{\Z_p} F^*D_Y\overset{\cong}\to \Q_p\otimes_{\Z_p}D_Y$$ which is 
compatible with $\nabla$. Here $F^*D_Y:=(F^*D)_Y$ and this is identified with the pullback $F_{D_X(Y)}^*(D_Y)$ of $D_Y$ where  $F_{D_X(Y)}: D_X(Y) \to D_X(Y)$ is the lifting of the Frobenius of $X$ induced by $F_Y:Y\to Y$.

We have
 $$D_X= \cO_X\otimes_{\cO_{D_X(Y)}} D_Y.$$

\end{para}

\begin{para}\label{lift}

If we are in the embedded situation and if $X=Y\otimes_{W(k)} k$ as a log scheme over $k$ (we call this a lifted situation), 
we have $\cO_{D_X(Y)}=\cO_Y$ and $D_X=  D_Y/pD_Y$. 

We have a lifted situation \'etale locally on $X$. 
\end{para}

\begin{para} Assume we are in the embedded situation \ref{XZ}. 

Following  \cite{O0} and \cite{O}, we define filtrations on $D_Y$ and on $F^*D_Y$ as follows. For $r\in \Z$, let
$$N_r(D)_Y:= D_Y\cap p^{-r}\Phi(F^*D_Y), \quad N^r(D)_Y=N_{-r}(D)_Y=D_Y\cap p^r\Phi(F^*D_Y),$$
$$M^r(D)_Y:= \{x\in F^*D_Y\;|\; p^{-r}\Phi(x)\in D_Y\}\subset F^*D_Y.$$
Thus  $N_{\bullet}$ is an increasing filtration on $D_Y$ and $N^{\bullet}$ (resp. $M^{\bullet}$) is a decreasing filtration on $D_Y$ (resp. $F^*D_Y$). 
Though Ogus uses $N^{\bullet}$, we mainly use $N_{\bullet}$ in this paper  because  the presentation of the isomorphism in \ref{rr} becomes simpler and 
the relation to the numbering in the theory of Frobenius gauges (\ref{4.1}) is better for $N_{\bullet}$. But $N^{\bullet}$ has a better relation with Deligne's  d\'ecal\'e operation  Dec$^r$ (\cite{De} 1.3, see \ref{Dec}). 

We have

\begin{sbpara}\label{ggll}
(1) For $r\gg 0$, $N_r(D)_Y=D_Y$.
For  $r\ll 0$, $N_r(D)_Y= p^{-r}\Phi(F^*D_Y)$ and hence $N_{r-1}(D)_Y= pN_r(D)_Y$.

(2) For $r\ll 0$, $M^r(D)_Y=F^*D_Y$. For $r\gg 0$, $M^{r+1}(D)_Y=pM^r(D)_Y$. 
\end{sbpara}
We have an isomorphism 

\begin{sbpara}\label{rr}
$$p^{-r}\Phi: M^r(D)_Y\overset{\cong}\to N_r(D)_Y\quad \text{for all} \;\;r\in \Z.$$
\end{sbpara}

For the $F$-crystal $\bf 1$ (\ref{dual}),  

\begin{sbpara}\label{1NM}  $N_r({\bf 1})$ is $\cO_Y$ if  $r\geq 0$ and is $p^{-r}\cO_Y$ if  $r\leq 0$. $M^r({\bf 1})_Y$ is $\cO_Y$ if $r\leq 0$ and $p^r\cO_Y$ if $r\geq 0$. 
\end{sbpara}
For Tate twists, we have 

\begin{sbpara}\label{rs}
$$N_r(D(s))_Y= N_{r+s}(D)_Y, \quad M^r(D(s))_Y= M^{r+s}(D)_Y. $$
\end{sbpara}

\end{para}

\begin{lem}\label{indep} (\cite{O0} 1.10.1.) Let $r\in \Z$. 

(1) The image of $N_r(D)_Y \to D_X$ is independent of the choice of an embedded situation $X\subset Y$. 

(2) The image of $M^r(D)_Y\to F^*D_X$ is independent of the choice of an embedded situation $X\subset Y$. 

\end{lem}

\begin{pf} 
We prove (1).   
Let $X\subset Y$ and let $X\subset Y'$ be the embedded situations,  let $Z:=Y\times Y'$ with the lifting $F_Z:=F_Y\times F_{Y'}$ of Frobenius, and consider the diagonal embedding $X\subset Z$.
 Since $\cO_{D_X(Y)}/p^n\cO_{D_X(Y)}\to \cO_{D_X(Z)}/p^n\cO_{D_X(Z)}$ and $\cO_{D_X(Y')}/p^n\cO_{D_X(Y')}\to \cO_{D_X(Z)}/p^n\cO_{D_X(Z)}$ are  flat for all $n\geq 0$ (this is because $\cO_{D_X(Z)}/p^n\cO_{D_X(Z)}$ is isomorphic locally to a PD polynomial ring over $\cO_{D_X(Y)}/p^n\cO_{D_X(Y)}$ and similarly for $Z\to Y'$), $D_Z \cap p^{-r}\Phi(F^*D_Z)$ is the module theoretic pullback of $D_Y \cap p^{-r}\Phi(F^*D_Y)$ and is the module theoretic pullback of $D_{Y'} \cap p^{-r}\Phi(F^*D_{Y'})$. This proves (1).

The proof of (2) is given in the similar way. 
\end{pf}

\begin{para}\label{Xr}
We will denote the independent image of $N_r(D)_Y\to D_X$ (resp. $M^r(D)_Y\to F^*D_X$) (\ref{indep}) by $N_r(D_X)$ (resp. $M^r(D_X)$). It is defined globally on $X$ and is a coherent subsheaf of $D_X$ (resp. $F^*D_X$). 

\end{para}

\begin{para}  By Ogus, the following two conditions (i) and (ii) are equivalent.

(i) For every $r\in \Z$, $N_r(D_X)$ is locally a direct summand of $D_X$.

(ii) For every $r\in \Z$, $M^r(D_X)$ is locally a direct summand of $F^*D_X$.

Ogus calls $D$ a uniform $F$-crystal when these equivalence conditions are satisfied.

The equivalence of (i) and (ii) is proved as follows.

\end{para}
\begin{prop}\label{adpro} The following conditions (i), (i)', (ii), (ii)', (ii$^{\star}$), (ii$^{\star}$)'  are equivalent.

(i) (resp. (i)') For every $r\in \Z$, the $\cO_X$-module 
$N_r(D_X)$ (resp. $M^r(D_X)$)  (\ref{Xr}) is locally a direct summand of $D_X$ (resp. $F^*D_X$).

(ii)  (resp. (ii)') \'Etale locally on $X$, for some lifted situation \ref{lift}, there is \'etale locally a direct sum decomposition $D_Y=\oplus_{r\in \Z}\; N_{\langle r\rangle}$ (resp. $F^*D_Y= \oplus_{r\in \Z} \; M^{\langle r\rangle}$) as an $\cO_Y$-module such that $$\Phi(F^*D_Y)=\oplus_{r\in \Z} \; p^rN_{\langle r\rangle}\quad (\text{resp.}\;\;D_Y= \oplus_{r\in \Z}  \; p^{-r}\Phi(M^{\langle r\rangle}))\quad \text{for all}\;\;r\in \Z.$$

The condition (ii$^{\star}$) (resp. (ii$^{\star}$)') is the one obtained from (ii) (resp. (ii)') by replacing ``for some lifted situation \ref{lift}, there is'' by ``for each lifted situation \ref{lift}, there is \'etale locally''.

\end{prop}

\begin{para}\label{pfad1} We prove \ref{adpro} in \ref{pfad1}--\ref{pfad8}. 

First, the equivalences (ii) $\Leftrightarrow$ (ii)' and   (ii$^{\star}$) $\Leftrightarrow$ (ii$^{\star}$)' 
follow from \ref{rr}. ($N_{\langle r\rangle}$ in (ii) and (ii$^{\star})$ and $M^{\langle r\rangle}$ in (ii)' and (ii$^{\star}$)' are connected by $N_{\langle r\rangle}=p^{-r}\Phi(M^{\langle r\rangle})$.) The implications (ii$^{\star}$) $\Rightarrow$ (ii), (ii$^{\star}$)'$\Rightarrow$ (ii)' are clear.

\end{para}

\begin{lem}\label{Drseq} 
(See \cite{O0} 1.12.4.) Let $X\subset Y$ be an embedded situation. For $r\in \Z$, we have the following exact sequences.

(1) $0 \to N_{r+1}(D)_Y\overset{p}\to N_r(D)_Y \to (N_r(D)_Y+pD_Y)/pD_Y\to 0.$

(2) $0 \to M^{r-1}(D)_Y\overset{p}\to M^r(D)_Y \to (M^r(D)_Y+pF^*D_Y)/pF^*D_Y\to 0.$

(3) $0 \to N_{r+1}(D)_Y/N_r(D)_Y\overset{p}\to N_r(D)_Y/N_{r-1}(D)_Y\to (N_r(D)_Y+pD_Y)/(N_{r-1}(D)_Y+pD_Y)\to 0.$  

(4) $0 \to M^{r-1}(D)_Y/M^r(D)_Y\overset{p}\to M^r(D)_Y/M^{r+1}(D)_Y\to (M^r(D)_Y+pF^*D_Y)/(M^{r+1}(D)_Y+pF^*D_Y)\to 0.$ 

\end{lem}
\begin{pf}
(1) and (2) are straightforwards. (3) follows from (1), and (4) follows from (2). 
\end{pf}

\begin{para}\label{pfad5} We prove the implication 
(i) $\Rightarrow$ (ii$^{\star}$) in  \ref{adpro}.

Consider a lifted situation $X\subset Y$. 
Note that for $r\in \Z$,    $(N_r(D)_Y+pD_Y)/(N_{r-1}(D)_Y+pD_Y)=N_r(D_X)/N_{r-1}(D_X)$ is locally free by our assumption. Locally, take a base  $(e_{r,i})_{1\leq i\leq  m(r)}$ of  $N_r(D_X)/N_{r-1}(D_X)$.  Locally, 
lift this base to a family of local sections $\tilde e_{r,i}$  of $N_r(D)_Y$. 
Then  $(\tilde e_{r,i})_{r\in \Z,1\leq i\leq m(r)}$ is a base of $D_Y$. 
For $r\in \Z$, let $N_{\langle r\rangle}$ be the $\cO_Y$-submodule of $D_Y$ generated by $\tilde e_{r,i}$ for $1\leq i\leq m(r)$. Then $D_Y$ is the direct sum of $N_{\langle r\rangle}$. 

For $r\in \Z$, let $$S_r= (\oplus_{s\leq r} \; N_{\langle r\rangle}) \oplus (\oplus_{s>r}\; p^{s-r}N_{\langle r\rangle}).$$  We prove $N_r(D)_Y=S_r$.
We have $$S_r\subset N_r(D)_Y, \quad S_r\subset S_{r+1}, \quad pS_{r+1}\subset S_r\;\;(r\in \Z).$$
Consider the commutative diagram 
$$\begin{matrix} 0&\to & S_{r+1}/S_r & \overset{p}\to & S_r/S_{r-1}& \to & (S_r+pD_Y)/(S_{r-1}+pD_Y) & \to 0\\
&&\downarrow &&\downarrow && \downarrow \\
0&\to & N_{r+1}(D)_Y/N_r(D)_Y & \overset{p}\to & N_r(D)_Y/N_{r-1}(D)_Y& \to & N_r(D_X)/N_{r-1}(D_X) & \to  0.\end{matrix}$$
The lower sequence is that of \ref{Drseq} (3) and hence is exact. The upper sequence is exact as is easily seen. The right vertical arrow is an isomorphism as is easily seen. 
Hence the middle vertical arrow is an isomorphism if and only if the left vertical arrow is exact. For $r\gg 0$, $S_r=N_r(D)_Y=D_Y$. Hence by the downward induction on $r$, we have that $S_r/S_{r-1}\overset{\cong}\to N_r(D)_Y/N_{r-1}(D)_Y$. By using the fact $S_r=N_r(D)_Y=D_Y$ for $r\gg 0$ and by using the downward  induction on $r$, we have $S_r=N_r(D)_Y$ for all $r$. 

For $r\ll 0$,  $S_r=\oplus_{s\in \Z} \; p^{s-r}N_{\langle s\rangle}$.    
For $r\ll 0$, we have $p^{-r}\Phi(F^*D_Y)= p^{-r}\Phi(M^r(D)_Y)= N_r(D)_Y=S_r$ and  hence we have $\Phi(F^*D_Y)= \oplus_{s\in \Z} \; p^sN_{ \langle s\rangle}$. Hence the condition (1$^{\star}$) is satisfied.

The proof of the implication (i)' $\Rightarrow$ (ii$^{\star}$)' is given in a similar way. 

\end{para}

\begin{para}\label{pfad8} We prove the implication (ii) $\Rightarrow $ (i) of \ref{adpro}. Assume we are in the lifted situation and assume $D_Y=\oplus_r N_{\langle r\rangle}$ as in (ii). Then $$D_X=\oplus_r \; N_{\langle r\rangle}/pN_{\langle r\rangle}
\supset N_r(D_X)=\oplus_{s\leq r} \; N_{\langle s\rangle}/pN_{\langle s\rangle}.$$
Hence the condition (i) is satisfied.

The proof of the implication  (ii)' $\Rightarrow $ (i) is given in a similar way. 
\end{para}

\begin{prop} The category of uniform $F$-crystals is stable under $\oplus$, $\otimes$ and the dual, and Tate twists.

\end{prop}

\begin{pf} The condition (ii)  in \ref{adpro} is stable under tensor products, the dual, and Tate twists. \end{pf}

\begin{para}\label{pbk} (See \cite{O0} 1.12.2,  1.12.3). Assume $D$ is uniform. Then for each $r\in \Z$, there is a unique crystal $N_r(D)$ (resp. $M^r(D)$) on the log crystalline site of $X$ over $W(k)$ which is locally free of finite rank endowed with a homomorphism $N_r(D)\to  D$ (resp. $M^r(D)\to F^*D$) satisfying the following condition. \'Etale locally on $X$, for each embedded situation, 
the induced homomorphism $N_r(D)_Y\to D_Y$ (resp. $M^r(D)_Y \to F^*D_Y$) induces an 
isomorphism $N_r(D)_Y\overset{\cong}\to D_Y\cap p^{-r}\Phi(F^*D_Y)$ (resp. $M^r(D)_Y\overset{\cong}\to \{x\in F^*D_Y\;|\; p^{-r}\Phi(x)\in D_Y$).

This  is seen as follows. Let $X\overset{\subset}\to Y$ be a lifted situation having the property stated in (ii). Then for $r\in \Z$, we have 
$$N_r(D)_Y= ( \oplus_{s\leq  r}  \; N_{\langle s\rangle}) \oplus (\oplus_{s>r}\; p^{s-r}N_{\langle s\rangle}).$$
Hence $N_r(D)_Y$ is locally free of finite rank. Since $N_r(D)_Y$ is stable under the connection in $D_Y$, $N_r(D)_Y$ determines an $F$-crystal on the log crystalline site of $X$ over $W(k)$ which is locally free of finite rank. This $N_r(D)$ has the desired property on $Y$ on this $Y$. 

Next we consider embedded situations. Let $D'$ be a crystal on the log crystalline site which is locally free and of finite type endowed with a homomorphism $D'\to D$. For an embedded situation, consider the 
following condition $(*)_Y$.

$(*)_Y$. The induced map $D'_Y\to D_Y$ induces an isomorphism $D'_Y\overset{\cong}\to D_Y\cap p^{-r}\Phi(F^*D_Y)$. 

We show that for embedded situations $X\subset Y$ and $X\subset Y'$, the conditions $(*)_Y$ and $(*)_{Y'}$ are equivalent. 
 Let $Z=Y\times Y'$, and let $X\subset Z$ be the diagonal embedding. Assume $(*)_Y$ is satisfied. Since $\cO_{D_X(Y)}/p^n\cO_{D_X(Y)}\to \cO_{D_X(Z)}/p^n\cO_{D_X(Z}$ are flat, $(*)_Z$ is satisfied. Since  $\cO_{D_X(Y')}/p^n\cO_{D_X(Y')}\to \cO_{D_X(Z)}/p^n\cO_{D_X(Z)}$ are faithfully flat, the condition $(*)_{Y'}$ is satisfied. 
 
 Now since $(*)_Y$ is satisfied locally  in some lifted situation as above, it is satisfied always.

\end{para}

\begin{prop}\label{NMdual}  Assume we are in a embedded situation. 

(1) Let $D$ and $D'$ be  uniform $F$-crystals $D$. Then we have  
$$N_r(D\otimes D')_Y=\sum_{i+j=r} N_i(D) \otimes N_j(D'),    \quad M^r(D\otimes D')_Y=\sum_{i+j=r} M^i(D)_Y \otimes M^j(D')_Y.$$

(2) Let $D$ be a uniform $F$-crystal $D$ and let $r\in \Z$. 

Then $N_r(D^*)_Y$ is the sheaf of linear maps $D_Y\to \cO_{D_X(Y)}$ which send $N_s(D)_Y$ to $N_{s+r}({\bf 1})_Y$ for all $s\in \Z$. That is, $N_r(D^*)_Y$ is the sheaf of linear maps $D_Y\to \cO_{D_X(Y)}$ which send $N_s(D)$ to $p^{-r-s}\cO_{D_X(Y)}$ for all $s\in \Z$. 

$M^r(D^*)_Y$ is the sheaf of  linear maps $F^*D_Y\to \cO_{D_X(Y)}$ which send $M^s(D)_Y$ to $M^{s+r}({\bf 1})_Y$ for all $s\in \Z$. That is, $M^r(D^*)_Y$ is the sheaf of linear maps $D_Y\to \cO_{D_X(Y)}$ which send $M^s(D)$ to $p^{r+s}\cO_{D_X(Y)}$ for all $s\in \Z$. 

\end{prop}

\begin{pf}
In the lifted situation, these are 
 seen by the conditions (ii$^{\star}$) and (ii$^{\star}$)' in \ref{adpro}. In fact, for the tensor product, if $D_Y=\oplus_r N_{\langle r\rangle}$ and $F^*D_Y= \oplus_r M^{\langle r\rangle}$ are the decompositions there and if we take a similar decompositions of $D'_Y$ and $F^*D'_Y$ by $N'_{\langle r\rangle}$ and $(M')^{\langle r\rangle}$, respectively,  $(D\otimes D')_Y$ has the decomposition whose $r$-th component is $\oplus_{i+j=r} N_{\langle i\rangle}\otimes N'_{\langle j\rangle}$ and $F^*(D\otimes D')_Y$ has the decomposition whose $r$-th component is  $\oplus_{i+j=r}M^{\langle i\rangle}\otimes (M')^{\langle j\rangle}$ which have the desired properties. $D^*_Y$ has the decomposition whose $r$-th component is the dual of $N_{\langle -r\rangle}$ and $F^*D^*_Y$ has the decomposition whose $r$-th component is the dual of $M^{\langle -r\rangle}$ having the desired properties. 
 
 In the embedded situation $X\subset Y$, take locally a lifted situation $X\subset Y'$, let $Z=Y \times Y'$, and let $X\subset Z$ be the diagonal embedding. By the fact $\cO_{D_X(Y')}/p^n\cO_{D_X(Y')}\to \cO_{D_X(Z)}/p^n\cO_{D_X(Z)}$ are flat, we have the above result in the case $Z$ is the $Y$ in the proposition, and then obtain the above result by the fact 
$\cO_{D_X(Y)}/p^n\cO_{D_X(Y)}\to \cO_{D_X(Z)}/p^n\cO_{D_X(Z)}$ are faithfully flat. 
\end{pf}

\begin{para}\label{Hodge0}

Assume $D$ is uniform. 

 In $F^*D_X$,  $M^r(D_X)$ (\ref{Xr}) is locally a direct summand of $F^*D_X$ and is stable under $\nabla$ of $F_X^*D_X$. Hence $M^r(D_X)$ descends to the $\cO_X$-submodule $$\fil^rD_X:=\text{Ker}(\nabla: M^r(D_X)\to M^r(D_X) \otimes_{\cO_X} \omega_X^1)$$ of $D_X$ which is locally a direct summand satisfying $M^r(D_X)=F^*\fil^rD_X$.

For uniform $F$-crystals $D$ and $D'$, we have 
$\fil^r(D\otimes D')_X=\sum_{s+t=r} \fil^sD_X\otimes_{\cO_X} \fil^tD'_X$. 
For the dual $D^*$ of a uniform $F$-crystal $D$, we have 
$\fil^rD^*_X={\cH}om_{\cO_X}(D_X/\fil^{-r+1}D_X,\cO_X)$. 

\end{para}

\begin{para}\label{MBO} By Mazur (\cite{Ma1}, \cite{Ma2}),  
Berthelot-Ogus \cite{BO} Theorem 8.26, and by Ogus  \cite{O} Theorem 7.3.1, 
 if $D$ comes from a proper log smooth scheme over $X$ as the higher direct image of the structure sheaf of the crystalline site, under mild assumptions, we have the following (1)--(3). 
 
(1) $D$ is uniform.

(2)  $(\fil^rD_X)_{r\in\Z}$ coincides with the geometric Hodge filtration on $D_X$ and $(M^r(D_X))_{r\in \Z}$ coincides with the geometric Hodge filtration on $F^*D_X$. .

(3)  $(N_r(D_X))_{r\in \Z}$  coincides (up to a suitable shift) with the geometric conjugate Hodge filtration of $D_X$.

\end{para}

\begin{para}\label{Dseq}

 Let $0\to D' \to D \to D'' \to 0$ be a sequence of uniform $F$-crystals on $X$. We say this sequence is exact if it satisfies the following conditions (i) and (ii). 
 
 \medskip
 
 (i) This sequence is  exact as a sequence of crystals. (Equivalently, this sequence is exact as a sequence of sheaves on the log crystalline site. Equivalently, locally on $X$, the sequence $0\to D'_Y\to D_Y \to D''_Y\to 0$ is exact in every embedded situation.) 
 
 (ii) The sequence $0\to \fil^rD'_X\to \fil^rD_X \to \fil^rD''_X \to 0$ is exact for every $r\in \Z$. 
 
 \medskip
 
We show that this sequence of uniform $F$-crystals is exact if and only if the following condition (iii) (resp. (iv)) is satisfied.

 (iii)  The sequence $0\to N_r(D') \to N_r(D) \to N_r(D'') \to 0$ is an exact sequence of crystals for every $r\in \Z$.
 
 (iv)  The sequence $0\to M^r(D') \to M^r(D) \to M^r(D'') \to 0$ is an exact sequence of crystals for every $r\in \Z$.

Note that the condition (iii) (resp. (iv))  is equivalent to the condition that in every local lifted situation, the sequence $0\to N_r(D')_Y\to N_r(D)_Y \to N_r(D'')_Y \to 0$ (resp. $0\to M^r(D')_Y \to M^r(D)_Y \to M^r(D'')_Y\to 0$) of sheaves is exact for every $r$. We prove the implication (iv) $\Rightarrow$ (i)+(ii). By taking $r\ll 0$, (iv) implies (i). By \ref{Drseq}  (2) in the lifted situation, (iv) implies  the exactness of the sequence $0\to F^*\fil^rD'_X\to F^*\fil^rD_X \to F^*\fil^rD''_X\to 0$ and hence the exactness of $0\to \fil^rD'_X\to \fil^rD_X\to \fil^rD''_X\to 0$. We prove (i)+(ii) $\Rightarrow$ (iv). Assume (i)+(ii). By \ref{Drseq} (4) and by induction on $r$, the sequence $0\to M^r(D')_Y/M^{r+1}(D')_Y\to M^r(D)_Y/M^{r+1}(D)_Y \to M^r(D'')_Y/M^{r+1}(D'')_Y\to 0$ is exact for every $r$. 
By induction on $n$, the sequence $0\to M^r(D')_Y/M^{r+n}(D')_Y \to M^r(D)_Y/M^{r+n}(D)_Y \to M^r(D'')_Y/M^{r+n}(D'')_Y\to 0$ is exact for every $r$ and every $n\geq 0$. By taking $\varprojlim_n$, we have (iv). The equivalence (iii) $\Leftrightarrow$ (iv) follows from  \ref{rr}.

\end{para}

\begin{para} We give an example of a sequence of $F$-crystals which satisfies the condition (i) in \ref{Dseq} but not the condition (ii). 

Let $D$ be a free module of rank $2$ over the structure sheaf of the log crystalline site of $X$ with base $(e_i)_{i=1, 2}$, and define $\Phi$ by $\Phi(1\otimes e_1)=e_1$ and $\Phi(1\otimes e_2)= e_2+p^{-n}e_1$ for some integer $n\geq 1$. 
We have a sequence of $F$-crystals  $0\to {\bf 1}\to D \to {\bf 1}\to 0$, where $\bf 1$ is as in \ref{dual}, the second arrow sends $1$ to $e_1$, the third arrow sends $e_2$ to $1$. This sequence satisfies the above condition (i), but does not satisfy the above condition (ii).  
This is because in a lifted situation, 
$N_0(D)_Y = \cO_Y e_1+p^n\cO_Y e_2$ and hence $N_0(D)_Y\to N_0({\bf1})_Y=\cO_Y$ is not surjective (the image is $p^n\cO_Y$ in ${\bf 1}_Y=\cO_Y$).  
\end{para}

\begin{para}\label{FV}  Example. 
 Let $G=(G_n)_n$ be a $p$-divisible group over $X$ ($G_n=\text{Ker}(p^n: G\to G)$) and let $D$ be the $F$-crystal on $X$ associated to $G$ in the covariant way. If $G=\Q_p/\Z_p$, $D={\bf 1}$. If $G=(\Q_p/\Z_p)(1)$, $D$ is the $1$-Tate twist ${\bf 1}(1)$ of $\bf 1$. Then $p\Phi$ comes from a homomorphism $F^*D \to D$ which is usually written as $F$,  and there is a homomorphism  $V: D\to F^*D$ such that $F V=p$ and $V F=p$ (that is, $\Phi V=1$ and $V\Phi=1$).  $D$ is uniform, $M^r(D)_Y=D_Y$ for $r<0$ and $M^r(D)=p^rV(D_Y)$ for $r\geq 0$ in the embedded situation, and $D_X/\fil^0D_X$ is identified with $\Lie(G)$.

\end{para}

\begin{para} We give  examples of  $F$-crystals which are not uniform. 

Consider 
the lifted situation \ref{lift}. We consider an $F$-crystal $D$ of rank $3$ defined as follows. Take $s_1, s_2\in \cO(Y)$ such that $ds_i \in p\omega^1_Y$ for $i=1,2$. Since $F_Y(\omega^1_Y)\subset p\omega^1_Y$, $1-F_Y: \omega^1_Y\to \omega^1_Y$ is an isomorphism. Let
$\omega_i= -(1-F_Y)^{-1}(p^{-1}ds_i)\in \omega^1_Y$ for $i=1,2$. 

Let $D_Y$ be the free $\cO_Y$-module of rank $3$ with base $(e_i)_{1\leq i\leq 3}$. Let $\Phi: F^*D_Y\to D_Y$ be the linear map corresponding to the matrix
$$\begin{pmatrix} p & 0& 0\\0&p& 0\\
s_1&s_2&p\end{pmatrix}$$
that is, $\Phi(1\otimes e_i) =pe_i+ s_ie_3$ for $i=1,2$ and $\Phi(1\otimes e_3)=pe_3$. Define the integrable connection $\nabla$ on $D_Y$ by
$$\nabla(e_i)= e_3\otimes \omega_i \;\;(i=1,2), \quad \nabla(e_3)=0.$$
Then $\Phi$ and $\nabla$ are compatible. We have $\Phi h=p^2$ and $h\Phi=p^2$ where $h:D_Y\to F^*D_Y$ is the linear map defined by $h(e_i)= p\otimes e_i-s_1\otimes e_3$ for $i=1,2$ and  $h(e_3)=p \otimes e_3$.
Hence we have an $F$-crystal $D=(D, \Phi)$.

We have $M^r(D)_Y=F^*D_Y$ for $r\leq 0$, $M^r(D)_Y=p^{r-2}h(D_Y)$ for $r\geq 2$, and $M^1(D)_Y$ is identified with  the $\cO_Y$-submodule $M$ of $\cO_Y^{\oplus 3}$ consisting of all elements $(a_i)_{1\leq i\leq 3}$ such that $a_1s_1+a_2s_2\in p\cO_Y$, where we identify $(a_i)_i$ with $\sum_i\; a_i\otimes e_i\in F^*D_Y$. We have that $M^2(D)_Y=h(D_Y)$ is generated by $p\otimes e_i-s_i\otimes e_3$ for $i=1,2$ and $p\otimes e_3$. 

(1) We give an example in which $M^1(D)_Y$ is not locally free. Consider the case $X=\Spec(k[T_1, T_2])$ and $Y=\text{Spf}(A)$ with $A$ the $p$-adic completion of $W(k)[T_1, T_2]$.  
Let $s_i=T_i^p$ for $i=1,2$. 

We have an $\cO_Y$-homomorphism $M\to \cO_Y/p\cO_Y$ which sends $(a_i)_i$ to $s_2^{-1}(a_1\bmod p)= -s_1^{-1}(a_2\bmod p)$. This homomorphisms sends $v=(s_2, -s_1,0)\in M$ to $1$. Hence if $M$ was locally free at the point $x$ of $X$ corresponding to the maximal ideal $(T_1,T_2)$ of $k[T_1, T_2]$, this $v$ should be a part of a base of $M_x$ and the $\cO_{Y,x}$-module $M_x/\cO_{Y,x}v$ should be free.  But we have an exact sequence of $\cO_{Y,x}$-modules 
$0 \to \cO_{Y,x}v\to M_x \to \cO_{Y,x}\oplus I\to 0$ where $I$ is the ideal of $\cO_{Y,x}$ generated by $s_1$ and $s_2$ and  the third arrow is $(a_i)_i\mapsto (a_3, b)$ with $b=p^{-1}(a_1s_1+a_2s_2)$. Hence $M^1(D)_Y$ is not locally free.

(2) We give an example in which $M^2(D_X)$ is not locally a direct summand of $D_X$. Let $X=\Spec(k[T])$ and $Y=\text{Spf}(A)$ where $A$ is the $p$-adic completion of $W(k)[T]$, and   take $s_1=s_2=T^p$. Then $M^2(D)_Y= h(D_Y)$ is generated by  $p\otimes e_i-T^p\otimes e_3$ for $i=1,2$ and $p\otimes e_3$ and hence the image $M^2(D_X)$ of $M^2(D)_Y\to F^*D_X$ is generated by $T^p\otimes e_3$ which is not locally a direct summand.

\end{para}

\section{De Rham complexes}\label{s:DR}

Let $X$ and $D=(D,\Phi)$ be as in \ref{XD}, We assume that $D$ is uniform. 
We consider  various de Rham complexes associated to $D$.

\begin{para}\label{Dbo} In general, assume we are given a filtered complex $C$, that is, a complex endowed with a decreasing filtration $({}^qC)_q$ by subcomplexes ${}^qC$.
Then define the complex $\overline{C}$ by
$$\overline{C}^q= \{x\in ({}^qC)^q\;|\; dx\in ({}^{q+1}C)^{q+1}\}$$ ($(-)^q$ denotes the part of degree $q$ of a complex) with the differentials  $d: \overline{C}^q\to \overline{C}^{q+1}$ induced by $d: ({}^qC)^q\to ({}^qC)^{q+1}$. 
This is the complex denoted by $\text{Dec}^0C$ in Deligne \cite{De} 1.3. 

By \cite{De} 1.3,  we have:

\begin{sbpara}\label{DBO} Assume that we are given a homomorphism of filtered complexes $C\to C'$ such that ${}^qC \to {}^qC'$ are 
quasi-isomorphisms for all $q$. Then the induced map $\overline{C}\to \overline{C'}$ is a quasi-isomorphism.

\end{sbpara}

\end{para}

\begin{para}\label{DRr} Assume we are in an embedded situation \ref{XZ}.

For $r\in \Z$, we define a subcomplex es $N_r\DR(D)_Y$ and $N^r\DR(D)_Y$ of the de Rham complex $$\DR(D)_Y:= [D_Y \overset{\nabla}\to D_Y\otimes_{\cO_Y}\omega^1_Y\overset{\nabla}\to D_Y \otimes_{\cO_Y}\omega^2_Y\overset{\nabla}\to \dots]$$  of $D_Y$ and a subcomplex $M^r\DR(D)_Y$ of the de Rham complex of $F^*D_Y$, as follows. We will regard these complexes as complexes of sheaves on the \'etale site of $X$.

Let $$N_r\DR(D)_Y:=\overline{C}\quad \text{with}\;\;{}^qC= \DR(N_{r-q}(D))_Y,$$ 
$$N^r\DR(D)_Y=N_{-r}\DR(D)_Y= \overline{C}\quad \text{with}\;\;{}^qC= \DR(N^{r+q}(D))_Y,$$
$$M^r\DR(D)_Y:= \overline{C}\quad \text{with} \;\; {}^qC=p^q\DR(M^{r-q}(D))_Y.$$

That is, the part of degree $q$ of $N_r\DR(D)_Y$ (resp. $N^r\DR(D)_Y$, resp. $M^r\DR(D)_Y$) is $$\{x\in N_{r-q}(D)_Y \otimes_{\cO_Y} \omega^q_Y\;|\; \nabla(x)\in N_{r-q-1}(D)_Y \otimes_{\cO_Y} \omega^{q+1}_Y\},$$ 
$$(\text{resp.}\quad \{x\in N^{r+q}(D)_Y \otimes_{\cO_Y} \omega^q_Y\;|\; \nabla(x)\in N^{r+q+1}(D)_Y \otimes_{\cO_Y} \omega^{q+1}_Y\},$$ 
$$\text{resp.}\quad \{x\in p^qM^{r-q}(D)_Y \otimes_{\cO_Y} \omega^q_Y\;|\; \nabla(x) \in p^{q+1}M^{r-q-1}(D)_Y \otimes_{\cO_Y} \omega^{q+1}_Y\}).$$

For $r\gg 0$, we have $N_r\DR(D)_Y= \DR(D)_Y$. On the other hand, if we define a subcomplex of $\DR(D)_Y$ by $$N_{(-\infty)}\DR(D)_Y=N^{(\infty)}\DR(D)_Y:= \overline{C}\quad \text{with ${}^qC= p^q\Phi(\DR(F^*D)_Y)$},$$  we have  $N_r\DR(D)_Y=p^{-r}N_{(-\infty)}\DR(D)_Y$ for $r\ll 0$.

For $r\ll 0$, we have  $M^r\DR(D)_Y= \overline{C}$ with ${}^qC= p^q\DR(F^*D)_Y$ and this complex is independent of $r\ll 0$. We denote this independent complex by $M^{-\infty}\DR(D)_Y$. 

For Tate twists, we have
$$N_r\DR(D(s))_Y= N_{r+s}\DR(D)_Y, \quad M^r\DR(D(s))_Y= M^{r+s}\DR(D)_Y.$$
\end{para}

\begin{lem}\label{DRrr} For $r\in \Z$, we  have an isomorphism $$p^{-r}\Phi:  M^r\DR(D)_Y\overset{\cong}\to N_r\DR(D)_Y.$$ 
We have also an isomorphism $\Phi: M^{-\infty}\DR(D)_Y\overset{\cong}\to N_{(-\infty)}\DR(D)_Y$. 
\end{lem}

\begin{pf} This follows from  \ref{rr} and \ref{DBO}.    \end{pf}

\begin{para}\label{Dec} For a filtered complex $C$, Deligne defines the subcomplex $\text{Dec}^rC$ of $C$ by $(\text{Dec}^rC)^q= \{x\in ({}^{q+r}C)^q\;|\; dx\in ({}^{q+r+1}C)^{q+1}\}$  (\cite{De} 1.3). That is, $\text{Dec}^rC=\overline{B}$ (\ref{Dbo}) where $B$ is the filtered complex defined by  ${}^qB= {}^{q+r}C$.

Regard $\DR(D)_Y$ as a filtered complex $C$  with the filtration 
 ${}^qC= \DR(N^q(D))_Y$. 
 
 Then  we have for all $r\in \Z$
$$N^r\DR(D)_Y=\text{Dec}^r\DR(D)_Y. $$

We have
$$N^{(\infty)}\DR(D)_Y= \text{Dec}^{(\infty)}\DR(D)_Y$$
where
$$\text{Dec}^{(\infty)}\DR(D)_Y:=  p^{-r}\text{Dec}^r\DR(D)_Y\quad \text{for}\; r\gg 0$$ 
which is independent of $r\gg 0$.

It seems that $M^r\DR(D)_Y$ is not nicely understood in terms of $\text{Dec}$. 
\end{para}

\begin{para}\label{dc1}

By \ref{DBO}, these local constructions in embedded situations  are glued by the method of simplicial hyper-covering as in \cite{HK} (2.18), and we have global objects $N_r\DR(D)$
 for $r\in \Z\cup\{(-\infty)\}$, and  $M^r\DR(D)$ 
  for $r\in \Z\cup\{-\infty\}$ in the derived category of sheaves of $W(k)$-modules on the \'etale site of $X$. 

\end{para}

\begin{para} Assume that we are in the lifted situation. Let
$$\eta: D_Y\to F_Y^*D_Y\;;\; x \mapsto 1 \otimes x$$
be the canonical embedding. 
Define a decreasing filtration $A^r(D)_Y$ on $D_Y$ by
$$A^r(D)_Y=\{x\in D_Y\;|\;\eta(x)  \in M^r(D)_Y\}.$$ 

This notation $A$ of the filtration  is taken from \cite{O0} and \cite{O}. 
\end{para}

The proofs of the following \ref{om4} and \ref{filD} will be given later. 

\begin{prop}\label{om4}

 Assume we are in the lifted situation.  
 Then we have the Griffiths transversality 
$$\nabla(A^r(D)_Y) \subset A^{r-1}(D)_Y \otimes_{\cO_Y} \omega^1_Y.$$ 
\end{prop}

Unlike $M^r(D)_Y$, $A^r(D)_Y$ need not correspond to a crystal. We need not have $\nabla(A^r(D)_Y)\subset A^r(D)_Y\otimes_{\cO_Y} \omega^1_Y$. 

\begin{prop}\label{filD} Assume that we are in the lifted situation. Let $r\in \Z$. Then the $\cO_Y$-module $A^r(D)_Y$ is locally free of finite rank, and the natural map $F_Y^*(A^r(D)_Y)\to M^r(D)_Y$ is an isomorphism.

\end{prop}

\begin{rem}\label{Onew} 
This Theorem \ref{filD} is proved in \cite{O0} and \cite{O} under the assumption that the lifting $F_Y$ of Frobenius is parallelizable in the sense in \cite{O0}, \cite{O}. 
Similarly, Proposition \ref{Hodge} and  Theorem \ref{DRrr2} are proved in \cite{O0} and \cite{O}
under the assumption that  $F_Y$ is parallelizable.
(In these results, Ogus does not assume that $D$ is uniform.)

 After the author sent this paper to Ogus, he found  new proofs of \ref{filD}, \ref{Hodge}, \ref{DRrr2} without assuming $F_Y$ is parallelizable which work in more general situation (relative situation without assuming $D$ is uniform). He also showed these results can be also deduced from the work of Shiho \cite{Sh}. See Appendix.

\end{rem}

\begin{cor}\label{Adual}  Assume we are in a lifted situation. 

(1) Let $D$ and $D'$ be  uniform $F$-crystals $D$. Then we have  
$$A^r(D\otimes D')_Y=\sum_{i+j=r} A^i(D) \otimes A^j(D').$$

(2) Let $D$ be a uniform $F$-crystal $D$ and let $r\in \Z$. 
Then
$A^r(D^*)_Y$ is the sheaf of  linear maps $D_Y\to \cO_Y$ which send $A^s(D)_Y$ to $A^{s+r}({\bf 1})_Y$ for all $s\in \Z$. That is, $A^r(D^*)_Y$ is the sheaf of linear maps $D_Y\to \cO_Y$ which send $A^s(D)$ to $p^{r+s}\cO_Y$ for all $s\in \Z$.

(3) Let $0\to D' \to D \to D'' \to 0$ be a sequence of uniform $F$-crystals satisfying the condition (i) in \ref{Dseq}. Then this sequence satisfies the condition (ii) in \ref{Dseq} if and only if the sequence 
$0\to A^r(D')_Y \to A^r(D)_Y \to A^r(D'')_Y \to 0$ is exact for every $r\in \Z$.

\end{cor}

\begin{pf}
(1) and (2) follow from \ref{filD} and \ref{NMdual}. 

(3) follows from \ref{filD} and the equivalence (ii) $\Leftrightarrow$ (ii)$_M$ in \ref{Dseq}. 
\end{pf}

By  \ref{om4} and \ref{filD}, we have

\begin{prop}\label{Hodge}  In the lifted situation, $\fil^rD_X$ (\ref{Hodge0}) coincides with the image of $A^r(D)_Y\to D_X$.

\end{prop}
See  \cite{O0} Theorem 2.2 and Remark \ref{Onew}. 
\begin{para}
We have  the Griffiths transversality
$$\nabla(\fil^rD_X)\subset \fil^{r-1}D_X\otimes_{\cO_X} \omega^1_X.$$ 
This is deduced from \ref{om4} and \ref{Hodge}. But this is proved in 
\cite{O0} Theorem 2.6  in the following way. To prove this, we may assume that we are in a lifted situation with parallelizable $F_Y$, and the result  \ref{Hodge} in the parallelizable case was enough for the proof of this.

\end{para}

\begin{para}\label{DRr2}  
Assume we are in a lifted situation.
By \ref{om4}, for $r\in \Z$, we have a subcomplex 
$$A^r\DR(D)_Y:=[A^r(D)_Y \overset{\nabla}\to A^{r-1}(D)_Y\otimes_{\cO_Y} \omega^1_Y \overset{\nabla}\to A^{r-2}(D)_Y\otimes_{\cO_Y} \omega^2_Y \to \dots]$$ of $\DR(D)_Y$.

For $r\ll 0$, we have $A^r\DR(D)_Y=\DR(D)_Y$. $(A^r\DR(D)_Y)_{r\in \Z}$ is a decreasing filtration on $\DR(D)_Y$. 

For Tate twists, we have $A^r\DR(D(s))_Y= A^{r+s}\DR(D)_Y$.

\end{para}
\begin{prop}\label{DRrr2} Assume we are in a lifted situation. For $r\in \Z$, the map $$\eta: A^r\DR(D)_Y\to M^r\DR(D)_Y\;;$$ $$\;\eta(x\otimes \omega)= \eta(x)\otimes F_Y(\omega)\quad (x\in A^{r-q}(D)_Y, \omega\in \omega^q_Y)\quad\text{in degree $q$}$$  is a  quasi-isomorphism.
In other words (\ref{DRrr}), we have a quasi-isomorphism
$$p^{-r}F\eta: A^r\DR(D)_Y \to N_r\DR(D)_Y\;;\;$$
$$ x\otimes \omega\mapsto p^{-r} F\eta(x)\otimes F_Y(\omega) \;\;(x\in A^{r-q}(D)_Y, \omega\in \omega^q_Y)\quad\text{in degree $q$}.$$ 
\end{prop}

See \cite{O} Theorem 7.2.1  and Remark \ref{Onew}. 

In \ref{om1}--\ref{endpf}, we prove \ref{om4}, \ref{filD} and \ref{DRrr2}. We first prove \ref{om4} (\ref{om1}--\ref{om5}).    
Then we prove \ref{filD} and \ref{DRrr2} together by induction on $r$.

\begin{para}\label{om1}  We recall the Cartier isomorphism. Let $\omega^q_{X,d=0}:=\Ker(\omega^q_X\overset{d}\to \omega^{q+1}_X)$. Then we have an isomorphism 
$$C: \omega^q_{X,d=0}/d\omega^{q-1}_X\overset{\cong}\to \omega^q_X$$
characterized by the property $C(t_0^pd\log(t_1)\wedge \dots d\log(t_q))= t_0d\log(t_1)\wedge \dots d\log(t_q)$ for all $t_0\in \cO_X$ and $t_i\in M_X$ for $1\leq i\leq n$ ($M_X$ denotes the log structure of $X$.) See \cite{KFI} Thm. 4.12 (1).

\end{para}

\begin{lem}\label{om2} Let $q\geq 0$, and let $(e_i)_i$ be a family of sections of $\omega^q_{X,d=0}$ whose classes in $\omega^q_{X=0}/d\omega^{q-1}_X$ form an $\cO_X^p$-base. Then $(e_i)_i$ is an $\cO_X$-base of $\omega^q_X$. 

\end{lem}

\begin{pf} 
\'Etale locally take $t_i\in M_X$ ($1\leq i\leq n$) such that $d\log(t_i)$ form a base of $\omega^1_X$.  For a subset $I=\{m(1), \dots, m(q)\}$ ($m(1)<\dots <m(q)$)  of $\{1,\dots, n\}$ of order $q$, let $\alpha_J=d\log(t_{m(1)})\wedge\dots \wedge d\log(t_{m(q)})$. Then $\alpha=(\alpha_J)_J$  is a base of $\omega^q_Y$. Furthermore, their classes form an $\cO_X^p$-base of $\omega^q_{X,d=0}/d\omega_X^{q-1}$ because $C(\alpha_J)=\alpha_J$. Let $e$ be a family of sections of  $\omega_{X,d=0}^q$ whose classes  form an $\cO_X^p$-base of $\omega^q_{X,d=0}/d\omega^{q-1}_X$. Then $e$ and $\alpha$ are connected by a section of $GL_s(\cO_X^p)$ where $s$ is the rank of $\omega^q_X$. This shows that  $e$ and the base $\alpha$ of $\omega^q_X$ are connected by a section of $GL_s(\cO_X)$, and hence $e$ is a base of $\omega^q_X$
\end{pf}

\begin{lem}\label{om3}  Assume we are in the lifted situation, let $(\mu_i)_i$ a base of $\omega^q_Y$, and let  $\nu_i:=p^{-1}F_Y(\mu_i)\in \omega^q_Y$. Then $(\nu_i)_i$ is a base of $\omega^q_Y$. 

\end{lem}

\begin{pf} We have $\nu_i\bmod p\in \omega^q_{X,d=0}$ and $C(\nu_i\bmod p)= \mu_i\bmod p$. Hence $\nu_i\bmod p$ form an $\cO_X^p$-base of $\omega^q_{X,d=0}/d\omega^{q-1}_X$. We apply \ref{om2}. 
\end{pf}

\begin{para}\label{om5} We prove \ref{om4}.
Let $(\mu_i)_i$ be a local base of $\omega^1_Y$ and let $x\in A^r(D)_Y$. Write $\nabla(x)=\sum_i y_i\otimes \mu_i$ with $y_i\in D_Y$. Then since $\Phi\eta(x)\in p^rD_Y$, we have $\nabla \Phi\eta(x)\in p^rD_Y\otimes_{\cO_Y} \omega^1_Y$. 
 On the other hand, $\nabla \Phi\eta(x)= \Phi\eta\nabla(x)= p\cdot \sum_i \Phi(y_i)\otimes \nu_i$, where $\nu_i=p^{-q}F_Y(\mu_i)$.  Since $\nu_i$ form a base of $\omega^1_Y$ (\ref{om3}), we have $y_i\in A^{r-1}(D)_Y$. 
\end{para}

Next we consider \ref{filD} and \ref{DRrr2}.

In the following, $S(\ref{filD}, r)$ denotes the statement of \ref{filD}  for one given  $r$, and  $S(\ref{DRrr2}, r)$ denotes the statement of \ref{DRrr2} for one given $r$.

\begin{para}\label{qi1}  Fix $r\in \Z$. 
Define complexes $A$, $B$ and a subcomplex $C$ of $B$ as follows.  
Let $$A:=A^r\DR(D)_Y/pA^{r-1}\DR(D)_Y, \quad B:=M^r\DR(D)_Y/pM^{r-1}\DR(D)_Y.$$ 
We have the canonical map $\eta: A\to B$. 

For $q\in \Z$, let $I^q$ be the image of $p^{q+1}M^{r-q-1}(D)_Y\otimes_{\cO_Y} \omega^q_Y$ in $B^q$. Let $C^q= I^q+\nabla(I^{q-1})$.

\end{para}

\begin{lem}\label{qi20}

(1) $\nabla: B^{q-1}\to B^q$ induces an isomorphism $I^{q-1}\overset{\cong}\to C^q/I^q$. 

(2) The complex $C$ is acyclic.

\end{lem}

\begin{pf}
For the proof of (4), the problem is the injectivity of $\nabla: I^{q-1}\to B^q/I^q$ but this is clear.

(2) follows from (1). 

\end{pf}

\begin{lem}\label{qi2}  Assume $S(\ref{filD}, r-q)$ and $S(\ref{filD}, r-q-1)$. 

(1) We have an isomorphism $$A^q\overset{\cong}\to B^q/C^q.$$ 

(2) The map $$A^{r-q}(D)_Y/pA^{r-q-1}(D)_Y \otimes_{\cO_X} (\omega^q_{X,d=0})^{(p)}\to B^q/I^q\;;\; x \otimes \omega\mapsto p^q\eta(x) \otimes \omega$$ is an isomorphism, where $(\omega^q_{X,d=0})^{(p)}$ denotes $\omega^q_{X,d=0}$ on which $a\in \cO_X$ acts as the multiplication by $a^p$.

(3) The isomorphism in (2)  induces an isomorphism  $$A^{r-q}(D)_Y/pA^{r-q-1}(D)_Y \otimes_{\cO_X} (d(\omega_X^{q-1}))^{(p)}\overset{\cong}\to C^q/I^q.$$

\end{lem}  

\begin{pf} We prove (2). We have 
$$B^q/I^q\overset{\cong}\leftarrow \{x\in M^{r-q}(D)_Y \otimes_{\cO_Y} \omega^q_Y\;|\;\nabla(x) \in pM^{r-q-1}(D)_Y \otimes_{\cO_Y} \omega^{q+1}_Y\}/pM^{r-q-1}(D)_Y \otimes_{\cO_Y}\omega^q_Y$$
$$\overset{\cong}\leftarrow  (M^{r-q}(D)_Y/pM^{r-q-1}(D)_Y)^{\nabla=0}\otimes_{\cO_X^p} \omega_{X,d=0}^q\overset{\cong}\leftarrow A^{r-q}(D)_Y/pA^{r-q-1}(D)_Y \otimes_{\cO_X} (\omega^q_{X,d=0})^{(p)}$$
where the first isomorphism is given by $p^q$, $\nabla=0$ denotes the kernel of $\nabla$, and the last arrow is $x\otimes \omega\mapsto \eta(x) \otimes \omega$ which is an isomorphism by  $S(\ref{filD}, r-q)$ and $S(\ref{filD}, r-q-1)$.

We prove (3).  We have a surjection $p^q: M^{r-q}(D)_Y/pM^{r-q}(D)\otimes_{\cO_X}\omega_X^{q-1}\to I^{q-1}$. By $S(\ref{filD}, r-q)$, $M^{r-q}(D)_Y/pM^{r-q}(D)\otimes_{\cO_X}\omega_X^{q-1}\cong A^{r-q}(D)_Y/pA^{r-q}(D) \otimes_{\cO_X} (\omega^{q-1}_X)^{(p)}$. Hence the image of $\nabla: I^{q-1}\to B^q/I^q\cong A^{r-q}(D)_Y/pA^{r-q-1}(D)_Y \otimes_{\cO_X} (\omega^q_{X,d=0})^{(p)}$ coincides with $A^{r-q}(D)_Y/pA^{r-q-1}(D)_Y \otimes_{\cO_X} (d(\omega_X^{q-1}))^{(p)}$. 

(1) follows from (2) and (3) by the fact that $p^{-q}F_Y: \omega^q_Y\to \omega^q_Y$ induces the inverse $\omega^q_X \overset{\cong}\to \omega^q_{X,d=0}/d\omega^{q-1}_X$ of the Cartier isomorphism.

\end{pf}

By \ref{qi20} and  \ref{qi2}, we have

\begin{cor}\label{qi4} (1) Assume  $S(\ref{filD},s)$ for all $s\leq r$. Then the map $A\to B$ is a quasi-isomorphism. 

(2) Assume $S(\ref{filD}, s)$ for all $s<r$. Then $\cH^q(C)=0$ for any $q\geq 1$, $A^q\overset{\cong}\to B^q/C^q$ for each $q\geq 1$, $\cH^q(A)\to \cH^q(B)$ is an isomorphism for each $q\geq 2$, and the map $\cH^1(A)\to \cH^1(B)$ is surjective.
\end{cor}

\begin{lem}\label{qi5} (1) Assume $S(\ref{filD}, s)$ for all $s\leq r$. Then we have $S(\ref{DRrr2}, r)$. 

(2) Assume $S(\ref{filD}, s)$ for all $s<r$. Then the map $\cH^q(A^r\DR(D)_Y)\to \cH^q(M^r\DR(D)_Y)$ is an isomorphism for $q\geq 2$ and surjective for $q=1$.
\end{lem}

\begin{pf} For $n\geq 0$, let $K^n= p^nA^{r-n}\DR(D)_Y$, $L^n= p^nM^{r-n}\DR(D)_Y$. Then \ref{qi5} follows from \ref{qi4} 
by the exact sequences 
$$0 \to K^n/K^{n+1}\to K^0/K^{n+1} \to K^0/K^n\to 0, \quad 0\to L^n/L^{n+1}\to L^0/L^{n+1}\to L^0/L^n\to 0$$
and by $$K^0=\varprojlim_n K^0/K^n, \quad L^0=\varprojlim_n L^0/L^n.$$
 \end{pf}

\begin{lem}\label{qi6} Let $$\begin{matrix} 0&\to &S& \to &T&\to& Q &\to& 0\\
&&\downarrow && \downarrow && \downarrow &&\\
0&\to &S' &\to& T'& \to& Q&\to &0\end{matrix}$$ be a commutative diagram of  exact sequence of complexes in an abelian category. 
Assume that the parts of degree $q$ of these complexes for $q<0$ are zero. Assume further the following (i)--(iii). 

\medskip
(i) The map $T\to T'$ is a  quasi-isomorphism. 

(ii) The map $H^q(S) \to H^q(S')$ is an isomorphism for $q\geq 2$, and it is surjective for  $q=1$. 

(iii) The map $H^0(Q)\to H^0(Q')$ is injective. 

\medskip

Then the maps $S\to S'$ and $Q\to Q'$ are quasi-isomorphisms.

\end{lem}

\begin{pf} Let $S''$ (resp. $T''$, resp. $Q''$) be the mapping fiber (the $-1$ shift of the mapping cone) of the map $S\to S'$ (resp. $T\to T'$, resp. $Q\to Q'$). By the assumption, $T''$ is acyclic. We prove that $S''$ and $Q''$ are also acyclic. We have $H^q(S'') \cong H^{q-1}(Q'')$ for all $q\in \Z$. The assumptions tell that  $H^q(S'')=0$ unless $q=0,1$ and $H^q(Q'')=0$ if $q\leq 0$. We have $H^1(S'')\cong H^0(Q'')=0$, $H^0(S'')  \cong H^{-1}(Q'')=0$.
\end{pf}

\begin{para}\label{pf2}  Since $A^r(D)_Y= D_Y$ for $r\ll 0$, $S(\ref{filD}, r)$ is true for $r\ll 0$. By \ref{qi5}, $S(\ref{DRrr2}, r)$ is true for $r\ll 0$.

We prove  $S(\ref{DRrr2}, r)$ assuming $S(\ref{filD}, s)$ and $S(\ref{DRrr2},s)$ for all $s<r$.

We apply \ref{qi6} by taking $S=A^r\DR(D)_Y$, $T=A^{r-1}\DR(D)_Y$, $Q=T/S$. $S'=M^r\DR(D)_Y$, $T'=M^{r-1}\DR(D)_Y$, $Q'=T'/S'$. 
The injectivity of $\cH^0(Q)\to \cH^0(Q')$ is by the definition of $A^r(D)_Y$. By \ref{qi6}, $S(\ref{DRrr2}, r)$ is true. 

\end{para}

\begin{para}\label{endpf} We prove $S(\ref{filD}, r)$ assuming $S(\ref{filD}, s)$ for all $s<r$ and $S(\ref{DRrr2}, s)$ for all $s\leq r$. 

We have an isomorphism (\ref{Drseq} (2)) $$M^r(D)_Y/pM^{r-1}(D)_Y\overset{\cong}\to (M^r(D)_Y+pF^*D_Y)/pF^*D_Y=M^r(D_X).$$ Hence (\ref{Hodge0}), the $\cO_X$-module 
$M^r(D)_Y/pM^{r-1}(D)_Y$ is generated by the kernel of the connection

\medskip

(1) $\nabla: M^r(D)_Y/pM^{r-1}(D)_Y\to M^r(D)_Y/pM^{r-1}(D)_Y
\otimes_{\cO_X} \omega^1_X.$

\medskip
\noindent
 It is sufficient to prove that this kernel coincides with $A^r(D)_Y/pA^{r-1}(D)_Y$.

Let the complexes $A, B, C$ be as in \ref{qi1}. Since $B^0/I^0$ is identified with the kernel of the above connection (1) and since $A^0=A^r(D)_Y/pA^{r-1}(D)_Y$ and $I^0=C^0$, it is sufficient to prove that the map $A^0\to B^0/C^0$ is an isomorphism. 
 By $S(\ref{DRrr2}, s)$ for all $s\leq r$, the map $A\to B$ is a quasi-isomorphism. Since $C$ is acyclic (\ref{qi20}),  the map $A\to B/C$ is a quasi-isomorphism. Furthermore we  have $A^q\overset{\cong}\to B^q/C^q$ for all $q\geq 1$. Hence
$A^0\to B^0/C^0$ is an isomorphism. 

\end{para}

\begin{prop}\label{triD} Let $0\to D' \to D \to D''\to 0$ be an exact sequence of uniform $F$-crystals (\ref{Dseq}). We have distinguished triangles
$$N_r\DR(D') \to N_r\DR(D)\to \DR^N_r(D'') \to N_r\DR(D')[1] \quad\text{for $r\in \Z\cup{(-\infty)}$},$$  $$M^r\DR(D') \to M^r\DR(D)\to M^r\DR(D'') \to M^r\DR(D')[1] \quad\text{for $r\in \Z\cup\{-\infty\}$}.$$

\end{prop}

\begin{pf} We have to be careful that the construction $({}^qC)_q\mapsto \overline{C}$ in \ref{Dbo}  is not an exact functor 

The case of $N_r\DR$ follows from that of $M^r\DR$ by \ref{DRrr}. Consider the embedded situation. Let $C$ be the mapping cone of $M^r\DR(D')_Y \to M^r\DR(D)_Y$. It is sufficient to prove that the map  $C\to M^r\DR(D'')_Y$ is a quasi-isomorphism. We may assume that this is a lifted situation. Then by \ref{DRrr2}, we have quasi-isomorphisms $A^r\DR(D') \to M^r\DR(D')$, $A^r\DR(D) \to M^r\DR(D)$, and $ A^r\DR(D'')\to M^r\DR(D'')$. Furthermore, the sequence  $0\to A^r\DR(D')_Y \to A^r\DR(D)_Y\to A^r\DR(D'')_Y\to 0$ is exact (\ref{Adual}  (3)). This proves that 
the map $C\to M^r\DR(D'')_Y$ is a quasi-isomorphism. 
\end{pf}

\begin{para}\label{minus} We consider de Rham complexes with minus log poles.

Let $X$ be a smooth scheme over $k$ endowed with a  normal crossing divisor $E$. We endow $X$ with the log structure associated $E$.

In the embedded situation, let $\cO_{D_X(Y)}(-\log)$ be the invertible ideal of $\cO_{D_X(Y)}$ generated \'etale locally by any  local section $f$  of the log structure $M$ of the PD-envelop $D_X(Y)$ such that the image of $f$ under $M\to \cO_X$ defines $E$. 
We have minus log pole versions of the modified de Rham complexes by
$${}^-\DR(D)_Y:=\cO_{D_X(Y)}(-\log)\otimes_{\cO_{D_X(Y)}} \DR(D)_Y,$$
$${}^-N_rDR(D)_Y:=  \cO_{D_X(Y)}(-\log)\otimes_{\cO_{D_X(Y)}} N_r\DR(D)_Y\quad\text{for $r\in \Z\cup\{(-\infty)\}$,}$$ $${}^-M^r\DR(D)_Y:= \cO_{D_X(Y)}(-\log)\otimes_{\cO_{D_X(Y)}} M^r\DR(D)_Y\quad\text{for $r\in \Z\cup\{-\infty\}$}.$$ 
 We have an isomorphism $p^{-r}\Phi: {}^-M^r\DR(D)_Y\overset{\cong}\to {}^-N_r\DR(D)_Y$.  
 
 If we regard ${}^-\DR(D)_Y$ as a filtered complex $C$ with the filtration defined by ${}^qC={}^-\DR(N^q(D))_Y$, we have
 ${}^-N_r\DR(D)_Y=\text{Dec}^{-r}({}^-\DR(D))$. Define $\text{Dec}^{(\infty)}({}^-DR(D)_Y)= p^{-r}\text{Dec}^r({}^-\DR(D)_Y)$ for $r\gg 0$ which is independent of $r\gg0$. Then ${}^-\N_{(-\infty)}\DR(D)_Y=\text{Dec}^{(\infty)}({}^-\DR(D)_Y)$.

 These local constructions are glued to global objects ${}^-\DR(D)$, ${}^-N_r\DR(D)$, ${}^-M^r\DR(D)$ of the derived category.

In the lifted situation, we have the minus log pole version $${}^-A^r\DR(D)_Y:=\cO_Y(-\log) \otimes_{\cO_Y} A^r\DR(D)_Y.$$ 
\end{para}

\begin{sbpara}\label{-fil} Assume that we are in the lifted situation.  Then the map
$\eta: {}^-A^r\DR(D)_Y \to {}^-M^r\DR(D)_Y$ is a quasi-isomorphism.

\end{sbpara}

\begin{pf} The method of the proof of \ref{DRrr2} works for this minus log version. Note that we have already proved \ref{filD}. The proof of \ref{filD}$\Rightarrow$ \ref{DRrr2} in \ref{qi5} (1) works also for the minus log pole version by using the following minus log version \ref{-Cart} of the Cartier isomorphism \ref{om1}. 
\end{pf}

\begin{sbpara}\label{-Cart} 
Let $${}^-\omega^q_X:= \Omega^q_X(-\log E)=\cO_Y(-E) \otimes_{\cO_X} \Omega^q_X(-\log E)$$ be the sheaf of differential forms with minus log poles. 

We have an isomorphism 
$$C: {}^-\omega^q_{X,d=0}/d({}^-\omega^{q-1}_X)\overset{\cong}\to {}^-\omega^q_X$$
characterized by the property $C(t_0^pd\log(t_1)\wedge \dots d\log(t_q))= t_0d\log(t_1)\wedge \dots d\log(t_q)$ for all $t_0\in \cO_X(-E)$ and $t_i\in M_Y$ for $1\leq i\leq n$ ($M_X$ denotes the log structure of $X$.)

\end{sbpara}
The proof of \ref{-Cart} is given in the same way as that of the usual (logarithmic) Cartier isomorphism \cite{KFI} Thm. 4.12 (1). 

In the same way as \ref{triD}, we have

\begin{sbpara}\label{-triD} For an exact sequence (\ref{Dseq}) of uniform $F$-crystals $0\to D'\to D \to D'' \to 0$, we have distinguished triangles
$${}^-N_r\DR(D') \to {}^-N_r\DR(D)\to {}^-N_r\DR(D'') \to {}^-N_r\DR(D')[1] \quad
\text{for $r\in \Z\cup\{(-\infty)\}$},$$ $${}^-M^r\DR(D') \to {}^-M^r\DR(D)\to {}^-M^r\DR(D'') \to {}^-M^r\DR(D')[1] \quad\text{for $r\in \Z\cup\{-\infty\}$}.$$

\end{sbpara}

\section{Syntomic complexes}\label{s:Syn}

Let $X$ be as in \ref{XD} and let $D$ be a uniform $F$-crystal on $X$. 
We define the  syntomic complex of $D$. 

Our method is similar to the method of Trihan in \cite{Tr}. In the case $D={\bf 1}(r)$ ($r\in \Z$), our definition of the syntomic complex coincides with that in the paper \cite{BMS}. 

\begin{para} Assume we are in the embedded situation \ref{XZ}. 

Define  the complex $\cS(D)_Y$ as the mapping fiber ($-1$ shift of the mapping cone) of the map $$1- \Phi\eta : N_0\DR(D)_Y \to N_{(-\infty)}\DR(D)_Y$$  where $1$ is the inclusion map $N_0\DR(D)_Y \to N_{(-\infty)}\DR(D)_Y$, and $\Phi\eta: N_0\DR(D)_Y \to N_{(-\infty)}\DR(D)_Y$ is the map induced from $\Q_p\otimes_{\Z_p} \DR(D)_Y\to \Q_p \otimes_{\Z_p} \DR(D)_Y$  whose degree $q$-part is $x\otimes \omega\mapsto \Phi(\eta(x)) \otimes F_Y(\omega)$, where $x\in D_Y$, $\omega\in \omega^q_Y$, $\eta(x)=1\otimes x \in F_{D_X(Y)}^*D_Y=F^*D_Y$.  We can say that $\cS(D)_Y$ is the mapping fiber of
$$1-\Phi \eta: \text{Dec}^0\DR(D)_Y \to \text{Dec}^{(\infty)}\DR(D)_Y$$
(see \ref{Dec}). 
The complex $\cS(D)_Y$ is isomorphic to the mapping fiber of $$1- \eta \Phi: M^0\DR(D)_Y \to M^{-\infty}\DR(D)_Y$$  where $1$ is the inclusion map.  In fact we have the isomorphism $\Phi$ from the latter to the former.

\end{para}

\begin{para} By \ref{DBO}, these local constructions in embedded situations  are glued and we have a global object
  $\cS(D)$ in the derived category of sheaves of abelian groups on the \'etale site of $X$. 

We have a distinguished triangle $$\cS(D) \to N_0\DR(D) \to N_{(-\infty)}\DR(D)\to \cS(D)[1].$$

\end{para}

 \begin{para}\label{triS} If $0\to D'\to D \to D'' \to 0$ is an exact sequence  (\ref{Dseq}) of uniform $F$-crystals, we have a distinguished triangle
 $$\cS(D') \to \cS(D)\to \cS(D'') \to \cS(D)[1].$$
 This follows from \ref{triD}.

 \end{para}

\begin{para}\label{cT1}

The following object $\cT(D)$ will be used in Section \ref{s:Tam} together with the above $\cS(D)$. 
By \ref{DBO}, we also have an object $\cT(D)$ of the derived category which is represented by the complex $$\cT(D)_Y:=\text{Dec}^{(\infty)}\DR(D)_Y/\text{Dec}^0\DR(D)_Y$$ $$= N_{(-\infty)}\DR(D)_Y/N_0\DR(D)_Y\overset{\cong}\leftarrow M^{-\infty}\DR(D)_Y/M^0\DR(D)_Y$$ in each embedded situation, where the last isomorphism is given by $\Phi$. In the lifted situation, by \ref{DRrr2}, $\cT(D)_Y$ is  quasi-isomorphic to the complex $(D_Y/A^{-q}(D)_Y \otimes_{\cO_Y} \omega^q_Y)_q$. 

 If $0\to D'\to D \to D'' \to 0$ is an exact sequence  of uniform $F$-crystals, by \ref{triD}, we have a distinguished triangle
 $$\cT(D') \to \cT(D)\to \cT(D'') \to \cT(D)[1].$$

\end{para}

\begin{para}\label{cT2}

In the case $D$ is the $F$-crystal associated to a $p$-divisible group $G=(G_n)_n$  in the covariant way (\ref{FV}), $\cT(D)$ is identified with $\Lie(G)$ put in degree $0$. This is seen as follows. 

In an embedded situation, the degree $q$ part of $M^{-\infty}\DR(D)_Y/M^0\DR(D)_Y$ is $0$ unless $q=0$, and the degree $0$ part is
$$\{x\in F^*D_Y\;|\; \nabla(x) \in pF^*D_Y\}/\{x\in V(D_Y)\;|\; \nabla(x) \in pF^*D_Y\}.$$ Hence $\cT(D)_Y$ is identified with $\cH^0(\cT(D))$, and there is a canonical homomorphism $\cH^0(\cT(D))\to (F^*D_X/V(D_X))^{\nabla=0}=\Lie(G)$ (\ref{FV}). We prove that the last arrow is an isomorphism. The problem is local and hence we may assume that we are in a lifted situation. In this situation, $\cT(D)_Y$ is quasi-isomorphic to $\DR(D)_Y/A^0DR(D)_Y$ (\ref{DRrr2}) whose degree $q$ part is zero unless $q=0$ and whose degree $0$ part is $D_Y/A^0(D_Y)\cong D_X/\fil^0D_X\cong \Lie(G)$.

\end{para}

\begin{para}\label{SD2}

In the lifted situation, by \ref{DRrr2}, $\cS(D)_Y$  is  quasi-isomorphic to the mapping fiber of  $$1-\Phi\eta: A^0\DR(D)_Y\to \DR(D)_Y.$$  We will call this mapping fiber also the syntomic complex but we will denote it  by ${\mathscr S}(D)_Y$ to distinguish it  from $\cS(D)_Y$. 

 In this situation, by \ref{DRrr2}, $\cT(D)_Y$ is quasi-isomorphic $\DR(D)_Y/A^0\DR(D)_Y$.
 \end{para}
 
 \begin{para} This is a preparation for \ref{H1sD}.

For uniform $F$-crystals $D$ and $D'$ on $X$, we have an abelian group $\Ext(D, D')$ classifying exact sequences $0\to D' \to E \to D \to 0$ of uniform $F$-crystals in the sense of \ref{Dseq} ( whose group law is given by Baer sums. 

Let the $F$-crystal ${\bf 1}$ be as in \ref{dual}.

We have $\Ext(D, D')= \Ext({\bf 1}, D^*\otimes D')$ where $D^*$ denotes the dual $F$-crystal of $D''$. 

\end{para}

\begin{thm}\label{H1sD} We have canonical isomorphisms $$\Hom({\bf 1}, D)\cong H^0(X, \cS(D)),\quad 
 \Ext({\bf 1}, D)\cong H^1(X, \cS(D)).$$

\end{thm}

\begin{pf} 

First we consider $\Hom({\bf 1}, D)$. Consider an embedded situation. We have a canonical embedding $\cH^0(\cS({\bf 1}))\subset \cO_Y$ and the section $1$ of $\cO_Y$ belongs to $H^0(X, \cS({\bf 1}))$. For a homomorphism ${\bf 1} \to D$, we have an element of $H^0(X, \cS(D))$ as the image of $1$ under the induced homomorphism $H^0(X, \cS({\bf 1}))\to H^0(X, \cS(D))$. To prove that $\Hom({\bf 1}, D) \to H^0(X, \cS(D))$ is an isomorphism, we may work \'etale locally on $X$ and hence we may assume that we are in the lifted situation and we can use the syntomic ccomplex ${\mathscr S}(D)_Y$ (\ref{SD2}). 
 Then $H^0(X, {\mathscr S}(D)_Y)$ is understood as the set of sections $e$ of $A^0(D_Y)$ such that $F\eta(e)=e$ and $\nabla(e)=0$. The map $\Hom({\bf 1}, D)\to H^0(X, {\mathscr S}(D)_Y)$ sends $h: {\bf 1}\to D$ to $e=h(1) \in H^0(X, {\mathscr S}(D)_Y)$ and it is clear that $h\mapsto e$ is bijective.

We next consider $\Ext({\bf 1}, D)$. 
For an exact sequence 
$0\to D \to E \to {\bf 1} \to 0$ of uniform $F$-crystals, we have a distinguished triangle $\cS(D) \to \cS(E) \to \cS({\bf 1}) \to \cS(D)[1]$ and hence a homomorphism $H^0(X, \cS({\bf 1}))\to H^1(X, \cS(D))$. By taking the image of $1\in H^0(X, \cS({\bf 1}))$ under this homomorphism, we have $\Ext({\bf 1}, D)\to H^1(X, \cS(D))$. 
By \'etale localization, we have a homomorphism $\cE xt({\bf 1}, D) \to \cH^1(\cS(D))$ of sheaves on the \'etale site of $X$. 

We show that it is enough to prove that this map $\cE xt({\bf 1}, D) \to \cH^1(\cS(D))$ is an isomorphism. 
We have a commutative diagram of exact sequences
$$\begin{matrix} 0 & \to & H^1(X, {\cH}om({\bf 1}, D))& \to &\Ext({\bf 1}, D) &\to &H^0(X, {\cE}xt({\bf 1}, D))& \to& H^2(X, {\cH}om({\bf 1}, D))\\
&&\downarrow&&\downarrow &&\downarrow&&\downarrow \\
0 & \to & H^1(X, \cH^0(\cS(D)))& \to &H^1(X, \cS(D)) &\to &H^0(X, \cH^1(\cS(D)))& \to& H^2(X, \cH^0(\cS(D)))\end{matrix}$$
Here the last arrow of  the upper row is as follows. Assume that we are given an element of $H^0(X, {\cE}xt({\bf 1}, D))$ which corresponds to 
an extension $0\to D \to E \to {\bf 1}\to 0$  given on an \'etale covering $U\to X$. \'Etale locally on $U\times_X U$, we have an isomorphism $\alpha: p_1^*E\overset{\cong}\to p_2^*E$ of extensions where $p_i$ is the $i$-th projection $U \times_X U\to U$. \'Etale locally on $U\times_X U \times_X U$,  if $q_i:U \times_X U \times_X U\to U$ denotes the $i$-th projection and $q_{ij}: U \times_X U \times_X U \to U\times_X U$ denotes the $(i,j)$-th projection, we have two  isomorphisms $q_{13}^*(\alpha)$ and $q_{23}^*(\alpha)\circ q_{12}^*(\alpha)$ from $q_1^*E $ to $q_3^*E$, and $q_{13}^*(\alpha)^{-1}\cdot q_{23}^*(\alpha)\circ q_{12}^*(\alpha)$ gives a global section of $\cH om({\bf 1}, D)$ on $U\times_X U \times_X U$ and this gives an element of $H^2(X, \cH om({\bf 1}, D))$. 

The first and the fourth vertical arrows in this diagram are isomorphisms. To prove that the second vertical arrow is an isomorphism, it is sufficient to prove that the third vertical arrow is an isomorphism.

 Hence we may assume that we are in the lifted situation and we can use the syntomic complex ${\mathscr S}(D)_Y$ (\ref{SD2}). 
 The map $\cE xt({\bf 1}, D) \to \cH^1({\mathscr S}(D)_Y)$ is described as follows. For an exact sequence $0\to D \to E \to {\bf 1}\to 0$, lift $1\in \cO_Y$ locally to $e\in A^0(E)_Y$. Then 
$e-F\eta e\in D_Y$, $\nabla w \in A^{-1}(D)_Y\otimes \omega^1_Y$, and we have $(e-F\eta e, \nabla e)\in \cH^1({\mathscr S}(D)_Y)$. 
 We have the converse map $\cH^1({\mathscr S}(D)_Y)\to \cE xt({\bf 1}, D)$ defined as follows. Let $(x, \omega)$ with $x\in D_Y$, $\omega\in A^{-1}(D)_Y \otimes \omega^1_Y$ be in $\cH^1({\mathscr S}(D)_Y)$. 
Let  $E=D\oplus A$. Extend $\Phi$ of $D$ to $E$ by $\Phi(0,1)= (x,1)$ and extend $\nabla$ of $D_Y$ to $E_Y$ by $\nabla(0,1)=\omega$. Then $\Phi$ and $\nabla$ are compatible and we have an exact sequence of uniform $F$-crystals $0 \to D \to E \to {\bf 1}\to 0$. It is seen easily that these maps between $\cE xt({\bf 1}, D)$ and $\cH^1({\mathscr S}(D)_Y)$ are the converses of each other. 
\end{pf}

\begin{para}\label{sminus}  We have the minus log pole versions of $\cS(D)$ and $\cT(D)$ as follows. 

Assume $X$ is a smooth scheme over $k$ endowed with a  normal crossing divisor $E$.

In the embedded situation, define the minus log pole version ${}^-\cS(D)_Y$ of $\cS(D)_Y$ as the mapping fiber of 
$$1- \Phi\eta : \text{Dec}^0({}^-\DR(D)_Y) \to \text{Dec}^{(\infty)}({}^-\DR(D)_Y).$$ 
The local constructions are glued to a global object 
${}^-\cS(D)$ in the derived category. We have a distinguished triangle 
  $${}^-\cS(D) \to{}^-N_0\DR(D) \to{}^-N_{(-\infty)}\DR(D)\to {}^-\cS(D)[1].$$

In the embedded situation, let $${}^-\cT(D)_Y=\text{Dec}^{(\infty)}({}^-\DR(D)_Y)/\text{Dec}^0({}^-\DR(D)_Y).$$ 
 The local constructions are glued to a global object 
 ${}^-\cT(D)$ of the derived category. 

  If $0\to D'\to D \to D'' \to 0$ is an exact sequence  (\ref{Dseq}) of uniform $F$-crystals, by \ref{-triD}, we have a distinguished triangle
 $${}^-\cS(D') \to {}^-\cS(D)\to {}^-\cS(D'') \to {}^-\cS(D)[1], \quad {}^-\cT(D') \to {}^-\cT(D)\to {}^-\cT(D'') \to {}^-\cT(D)[1].$$

 In the lifted situation, by \ref{-fil}, ${}^-\cS(D)_Y$  is quasi-isomorphic to the mapping fiber ${}^-{\mathscr S}(D)_Y$ of 
 $$1-\Phi \eta : {}^-A^0\DR(D)_Y\to {}^-\DR(D)_Y$$
 and ${}^-\cT(D)_Y$ is quasi-isomorphic to the complex $(D_Y/A^{-q}(D)_Y \otimes_{\cO_Y} {}^-\omega^q_Y)_q$.

\end{para}

\section{On Tamagawa number  conjecture}\label{s:Tam}
 
 We give a class number formula  \ref{L(1)p}, which is an analogue in characteristic $p$ of Tamagawa number conjecture.

  \begin{para} Let $k$ be a finite field of characteristic $p$, and let $X$ be a proper smooth scheme over $k$ endowed with a simple normal crossing divisor $E$. (Here we assume $E$ is simple normal crossing, not only normal crossing, because in the result of Faltings which we use in \ref{L(u)}, he assumes it.) 
   Let $U:= X\smallsetminus E$. Let $\bar k$ be an algebraic closure of $k$.
   
   Let $D=(D,\Phi)$ be a uniform $F$-crystal on $X$.

   For the comparison, let $\ell$ be a prime number which is different from $p$, and let $T$ be a smooth $\Z_{\ell}$-sheaf on the \'etale site of $U$. 
  \end{para}

\begin{para}\label{L(u)}
We review the $L$-functions of $T$ and $D$.

(1) The $L$-function $L_U(T, u)$ of $T$ (defined as the Euler product over closed points of $U$) is equal to $\text{det}(1-\varphi_k u\;|\; R\Gamma_{et,c}(U\otimes_k \bar k,T_{\Q_{\ell}}))^{-1}$. Here $R\Gamma_{et,c}$ denotes the compact support \'etale cohomology and $\varphi_k$ is the map induced by the $\sharp(k)$-th power map on $X$.

(2) The $L$-function of $L_U(D,u)$ of $D$ (defined as the Euler product over closed points of $U$) is equal to $\text{det}(1-\varphi_k u\;|\; R\Gamma(X,{}^{-}\DR(D)_{\Q_p}))^{-1}$. Here $\varphi_k=\Phi^d$ with $d=[k:\F_p]$, and $R\Gamma(X, \;)$ is \'etale cohomology. Recall that ${}^-\DR(D)$ denotes the de Rham complex of $D$ with minus log poles along $E$.

For (1), see  \cite{SGA5}. For (2),  see \cite{Fa}.
\end{para}

We consider the $L$-value at $u=1$.

\begin{prop}\label{L(1)}

Assume that $1-\varphi_k : H^m_{et,c}(U\otimes_k \bar k, T_{\Q_{\ell}}) \to H^m_{et,c}(U\otimes_k \bar k, T_{\Q_{\ell}})$ are  isomorphisms for all $m$. Then $H^m_{et,c}(U, T)$ are finite for all $m$ and
$$L_U(T, 1) \equiv \prod_m \sharp(H^m_{et,c}(U, T))^{(-1)^m} \; \bmod \Z_{\ell}^\times.$$

\end{prop}

\begin{prop} \label{L(1)p}

Assume that $1-\Phi : H^m(X, {}^-\DR(D)_{\Q_p})\to H^m(X, {}^-\DR(D))_{\Q_p}$ are isomorphisms for all $m$. Then $H^m(X,{}^- \cS(D))$ and $H^m(X, ^-\cT(D))$ are finite for all $m$ and 
$$L_U(D, 1) \equiv \prod_m \sharp(H^m(X,{}^- \cS(D)))^{(-1)^m} \sharp(H^m(X, ^-\cT(D))^{(-1)^m}\; \bmod W(k)^\times.$$ 

Here ${}^-\cS(D)$ and ${}^-\cT(D)$ are as in \ref{sminus}.

\end{prop}

\begin{para}

This generalizes  the formula of Milne \cite{Mi} who treated the case where $D={\bf 1}(r)$ (\ref{dual})  for some $r\in \Z$.

\end{para}

 \ref{L(1)} and \ref{L(1)p} are deduced from \ref{L(u)} (1) and (2), respectively, as follows. 

\begin{para}\label{exact}

(1) We have an exact sequence  
$$\dots \to H^m_{et,c}(U, T)\to H^m_{et,c}(U\otimes_k \bar k, T)\overset{1-\varphi_k}\to H^m_{et,c}(U\otimes_k \bar k, T)\to H^{m+1}_{et,c}(U, T)\to \dots.$$

(2) Let $\cP= \text{Dec}^0({}^-\DR(D))$ and $\cQ=\text{Dec}^{(\infty)}({}^-\DR(D))$. We have an exact sequences  
$$\dots \to H^m(X, {}^-\cS(D))\to H^m(X, \cP)\overset{1-\Phi}\to H^m(X, \cQ)$$ $$\to H^{m+1}(X, {}^-\cS(D))\to \dots,$$
$$\dots \to H^m(X, \cP)\overset{1}\to H^m(X, \cQ)\to H^m(X, {}^-\cT(D))$$ $$\to H^{m+1}(X, \cP) \to \dots.$$

\end{para}

We give elementary lemmas.

\begin{lem}\label{det}
For a finitely generated $\Z_{\ell}$-module $H$ and for a $\Z_{\ell}$-homomorphism $h:H\to H$ which induces an isomorphism $H_{\Q_{\ell}}\overset{\cong}\to H_{\Q_{\ell}}$, we have 
$$\text{det}(h)\equiv  \sharp(\text{Coker}(h: H\to H))\cdot \sharp(Ker(h:H \to H))^{-1}\; \bmod \Z_{\ell}^\times.$$ 
\end{lem}

\begin{lem}\label{det2} Let $H$ and $H'$ be  finitely generated $W(k)$-module, let $\alpha: H'\to H$ be a $W(k)$-linear map such that $\alpha: H'_{\Q_p}\to H_{\Q_p}$ is an isomorphism, and  let $\beta: H'\to H$ be a Frobenius-linear map such that $\alpha-\beta: H'_{\Q_p} \to H_{\Q_p}$ is an isomorphism. 
Consider the $W(k)_{\Q_p}$-linear map $g:=(\beta\circ \alpha^{-1})^d: H_{\Q_p}\to H_{\Q_p}$ where $d=[k:\F_p]$. Then $1-g: H_{\Q_p}\to H_{\Q_p}$ is an isomorphism and we have 
$$\text{det}_{W(k)_{\Q_p}}(1-g)\equiv  \sharp(\text{Coker}(\alpha-\beta: H'\to H))
\cdot \sharp(Ker(\alpha-\beta: H'\to H))^{-1}$$ $$
\cdot \sharp(\text{Coker}(\alpha: H'\to H))^{-1}\cdot \sharp(\text{Ker}(\alpha: H'\to H))\; \bmod W(k)^\times.$$

\end{lem}

\begin{para} The zeta value formula \ref{L(1)} follows from the 
 presentation of $L_U(T, u)$ in \ref{L(u)} (1) and the exact sequence \ref{exact} (1),  by using  \ref{det} which we apply by taking $H=H^m_{et,c}(U\otimes_k \bar k, T)$ and 
$h=1-\varphi_k$. 

 The zeta value formula \ref{L(1)p} follows from 
 the presentation of $L_U(D, u)$ in \ref{L(u)} (2) and the exact sequence \ref{exact} (2),  by using \ref{det2} which we apply by taking $H'=H^m(X, \cP)$ and $H=H^m(X, \cQ)$, $\alpha$ (resp. $\beta$) the map indued from the inclusion map (resp. $\eta \Phi$) of de Rham complexes (so $g=\varphi_k$).  
\end{para}

\section{Complements}\label{s:Cpl}

The following \ref{4.1}--\ref{413} is about relations to the theory on Frobenius gauges \cite{FJ} of Fontaine and Jannsen. 
  \ref{4.3} is about Tamagawa number conjecture in characteristic $p$.

\begin{para}\label{4.1} Let $X$ be as in \ref{XD}. For simplicity, we assume here the log structure of $X$ is trivial. 

A Frobenius gauge on $X$ is a family $(D_r)_{r\in \Z}$ of crystals $D_r$  on the crystalline site of $X$ over $W(k)$ which are locally free of finite type (torsion crystals are also considered in the theory of Fontaine-Jannsen, but we consider only these $D_r$ for simplicity) endowed with homomorphisms 
$$f: D_r\to D_{r+1}, \quad v: D_r\to D_{r-1}\quad \text{for each} \;\; r\in \Z$$ satisfying
(1) $fv=p, vf=p$, (2) $f$ is an isomorphism $r\gg 0$,  (3) $v$ is an isomorphism for $r\ll 0$, and endowed with an isomorphism $$F^*D_{\infty}\overset{\cong}\to D_{-\infty}$$
where $$D_\infty= \underset{f}\varinjlim D_r, \quad D_{-\infty}= \underset{v}\varinjlim D_r.$$

 Let $D$ be a uniform $F$-crystal on $X$. Then the crystals $D_r:=N_r(D)$ ($r \in \Z$) in \ref{adpro} (iii) associated to $D$ form a Frobenius gauge on $X$ as follows: $f: D_r\to D_{r+1}$ is the homomorphism which gives the inclusion maps $N_r(D)_Y\to N_{r+1}(D)_Y$ in local embedded situations, $v:D_r\to D_{r-1}$ is the homomorphism which gives $p:N_r(D)_Y\to N_{r-1}(D)_Y$ in local embedded situations, and the isomorphism  $F^*D_{\infty} \overset{\cong}\to D_{-\infty}$ is obtained from 
 $D_{\infty}=D$ and $D_{-\infty}\cong F^*D$ (\ref{ggll} (1)).

 \end{para}
 
 \begin{para}\label{413}
 
As is illustrated in Introduction of \cite{FJ},   Fontaine and Jannsen have a generalization of Dieudenn\'e theory to Frobenius gauges. For $n\geq 1$, there is a functor 
 $\cF_n$ from the category of Frobenius gauges to the category of sheaves of $\Z/p^n\Z$-modules on the syntomic site (\cite{FM}) of $X$ which is related to the classical Diedonn\'e theory as follows. 
 
 Let $G=(G_n)_n$ be a $p$-divisible group over $X$ and let $D$ be the $F$-crystal on $X$ associated to $G$ in the covariant way (\ref{FV}). 
  Let $(D_r)_r$ be the associated Frobenius gauge. Then $$\cF_n((D_r)_r)= G_n.$$

For a general uniform $F$-crystal $D$ and the associated Frobenius gauge $(D_r)_r$,  we expect that for $\cF_n(D):= \cF_n((D_r)_r)$,  there is a canonical isomorphism 
 $$R\epsilon_* (\cF_n(D)) = \cS(D)\otimes_{\Z}^L \Z/p^n\Z$$
where $\epsilon$ is the morphism from the syntomic site  of $X$  to the \'etale site of $X$. 

 \end{para}

 \begin{para}\label{412} There is a definition of the 
 syntomic complex of a Frobenius gage $(D_r)_r$ which is compatible with the definition of the syntomic complex  in Section 4. In an embedded situation,  let $\cS((D_r)_r)_Y$ be the mapping fiber of $$1-\Phi \eta: \text{Dec}^0\DR(D_{\infty})_Y\to \text{Dec}^{(\infty)}\DR(D_{\infty})_Y$$ where the notation is as follows. We regard $D_{r,Y}\subset D_{r+1, Y}$ via $f$.
 Hence $D_{\infty,Y}=D_{r,Y}$ if $r\gg 0$. We identify $D_{-\infty, Y}$ with $p^rD_{r,Y}$ for $r\ll 0$ which is independent of $r\ll 0$. 
Regard $\DR(D_{\infty})_Y$ as a filtered complex $C$ with the filtration defined by ${}^qC= \DR(D_{-q})_Y$. Let $\text{Dec}$ be as in \ref{Dec}, and define $\text{Dec}^{(\infty)}\DR(D_{\infty})_Y:= p^{-r}\text{Dec}^r\DR(D_{\infty})_Y$ for $r\gg 0$ which is independent of $r\gg 0$. In the above, 
 $1$ in $1-\Phi\eta$ is the inclusion map. The map $\Phi\eta$ is defined by the composition $D_{\infty,Y} \overset{\eta}\to F^*D_{\infty,Y} \overset{\Phi}\to D_{-\infty,Y}$. 
 These local construction are glued to a global object $\cS((D_r)_r)$ in the derived category. 
 
  This definition of $\cS((D_r)_r)$ is given in the preprint \cite{K0} which is referred to in the book \cite{O}. 
  
  However, if we define the notion of exact sequence $0\to (D'_r)_r \to (D_r)_r \to (D''_r)_r\to 0$ of Frobenius gages by the exactness of $0\to D'_r \to D_r\to D''_r \to 0$ for all $r$ (this seems to be the natural definition of exactness), we may not have a distinguished triangle $\cS((D'_r)_r) \to \cS((D_r)_r) \to \cS((D''_r)_r) \to \cS((D'_r))[1]$ (the trouble is that the operator $\text{Dec}^r$ is not exact) and hence it may not be possible to have a good relation like Theorem \ref{H1sD} between Frobenius gauges and the cohomology of the associated syntomic complexes. This is the reason why the author discussed syntomic complexes of uniform $F$-crystals as in this paper not starting from Frobenius gauges. 
 
 \end{para}

\begin{para}\label{4.3}  We presented class number formulas \ref{L(1)} and \ref{L(1)p} as analogues of the Tamagawa number conjecture in \cite{BK}. But there are more things to be considered.

Assume $X$ is of dimension one and connected. 

(1) Our formulation is not very similar to that in \cite{BK}. In \cite{BK}, local measures and global Tamagawa measure were considered. In our situation, these measures  correspond to measures on $H^1(O_v, {}^-\cS(D))$ for closed points $v$ of $X$ where $O_v$ denotes the completion of the local ring $\cO_{X,v}$, and to the Tamagawa measure on the product $\prod_v H^1(O_v, {}^-\cS(D))$ where $v$ ranges over all closed points of $X$. These are not considered in this paper.

(2) In this paper, we considered  $F$-crystals which have log poles at points of bad reduction. Motives which have semi-stable reduction everywhere give such $F$-crystals. But this paper does not treat motives which need not have semi-stable reduction.

(3) In this paper, we considered only the situation where the zeta function does not have zero or pole and regulator is not involved. 

In the paper \cite{KT}, we treated any bad reduction of an abelican variety and treated regulator (height pairing). In the paper \cite{BT}, Brinon and Trihan allows zero and poles of zeta functions. 
\end{para}

\section{Appendix. Remarks on $F$-spans, by Arthur Ogus}\label{s:App}

K. Kato has pointed out 
that the hypothesis 
of ``parallelizability'' for the lifting $F$ of $F_X$
in \cite{O0} Theorem 2.2 is superfluous.
This note  is my attempt to understand his comment.

Let $X/S$ be a smooth morphism of schemes
in characteristic $p$, and let $F_{X/S} \colon X \to X'$
be the relative Frobenius morphism.
Suppose we are given a lifting
$F \colon Z \to Z'$ of $F_{X/S}$,
where $Z/T$ and $Z'/T$ are formally smooth formal liftings
of $X/S$ and $X'/S$ respectively and $T$ is $p$-torsion free.

Our main goal is the following result.

\begin{thm}\label{desce.t}
  Let $(E',\nabla')$ be a $p$-torsion free coherent sheaf
  with integrable
  and quasi-nilpotent connection on $Z'/T$.
  Suppose that $E$ is a  submodule of $\tE:=F^*E'$
  which is invariant under the induced connection
  $\tilde \nabla $ on $\tilde E$.    Let $\eta \colon E' \to F_*F^*E'$ be
  the adjunction map.  Then
  the natural map
  $F^*(\eta^{-1} E) \to E$ is an isomorphism.
\end{thm}

Since $F$ is faithfully flat, 
the most natural way to attack this problem would be
to show  that $E$ is necessarily invariant under the
descent data for $F^*E'$. In fact,
as we explain later, the result can be deduced from
Shiho's Theorem 3.1 in \cite{Sh} 
which uses related, but different, descent data.  
Here  we follow a different approach,
based on Cartier descent, working in the context
of F-spans, which in fact is where the above question found its origin.

We recall some notions from \cite{O0}.
Let $F \colon Z\to Z'$ be a lifting
of $F_{X/S}$ as above. Then a
\textit{$F$-span on $X/S$} is given by a pair of
$p$-torsion free
$\cO_Z$-modules with (nilpotent)
integrable connection $(E',\nabla')$ on $Z'/W$ 
and $(E,\nabla) $ on $Z/W$ together with an injective
homomorphism
$$\tilde \Phi \colon F^*(E',\nabla') \to (E,\nabla).$$
Typically one assumes that $E'$ and $E$ are locally
free of finite rank and that the image of $\tPhi$
contains $p^nE$ for some $n > 0$.  Here we assume
only the latter.

We have natural maps
$$\eta_F \colon E' \to F_*F^*E'$$
$$\Phi_F \colon E' \to E := \tilde \Phi_F \circ \eta_F.$$
We define:
\begin{eqnarray*}
M^i_FF^*E'& := &\tilde \Phi_F^{-1} (p^i E)\cr
A^i_F E'& := & \Phi_F^{-1} (p^i E) = \eta_F^{-1} (M_F^i F^*E') \cr
  M^{[i]}_FF^*E'& := &\sum_{j} p^{[j]} M^{i-j}_F F^*E' \cr
  A^{[i]}_F E'& := & \sum_{j} p^{[j]} A^{i-j}_FE' 
\end{eqnarray*}

We shall  show that the filtration $A^{[\cdot]}_F$, unlike
the filtration $A^\cdot_F$, 
 is independent of the lifting
$F$, allowing a conceptual simplification of some of the constructions of \cite{O}.

The map $p^{-i}\tilde \Phi_F$ induces a morphism
$M^i_FF^*E' \to E$; we denote by $N^{-i}_FE$ its image.
Thus we find an isomorphism
$$\tilde \Phi_{i} \colon M^{i}_F F^*E' \to N_F^{-i} E.$$
If there is no danger of confusion, we may
drop the subscript indicating the choice of the lifting $F$.
We write $E_X$ for $E/pE$, and similarly for $E'$ and $\tE$.

\begin{prop}\label{agft.p}
  The filtrations $M^{\cdot} \tE$, $M^{[\cx]}\tE$
  are stable under $\tnabla$, annd the filtration
  $N^\cx E$ is stable under $\nabla$.  
  The filtrations $A^\cx E'$ and $A^{[\cx]}E'$
  satisfy Griffiths transversality. The  isomorphisms
  $\tilde \Phi_i \colon M^i \tE\to  N^{-i} E$
  are compatible with the connections,
  and induce isomorphisms:
  \begin{eqnarray*}
(\gr^i_M \tE ,\tnabla ) & \to & (N^{-i}E_X,\nabla) \cr
                                (\gr^i_M \tE_X,\tnabla)  &\to& (\gr^{-i}_N E_X, \nabla)
\end{eqnarray*}                                                                 
\end{prop}
\begin{pf}
  Omitted.  See \cite{O}.
\end{pf}

\begin{thm}\label{fspan.t}
  With the notation above, and for every $i$, the following statements hold:
  \begin{enumerate}
 \item  The natural map
   $F^*(A^iE' )  \to M^i\tE$ is an isomorphism.

  \item The natural maps
    $F^*(A^iE'_X )\to M^i\tE_X$
and 
$A^iE'_X \to (M^iE_X)^\nabla$
are isomorphisms.

 \item  The natural maps
   $F^*\gr^i_A E'_X \to \gr^i_M \tE_X$ and
   $\gr^i_A E'_X \to (\gr^i_M \tE_X)^\tnabla$ are isomorphisms.
   \item The natural map
$F^*\gr^i_A E' \to \gr^i_M\tE$
is  an isomorphism.
  \end{enumerate}
\end{thm}
\begin{pf}
  We prove these statements together by induction on $i$.
  Suppose $i = 0$.  Then statement (1) is true by definition,
  and it follows that $F^*E'_X \cong \tE_X$.  Cartier descent implies that
  $E'_X = (F^*E'_X)^\tnabla$.
  This proves (2).

  Since $A^1E' = \eta^{-1}(M^1\tE)$ and
  contains $pE'$, it follows that
  $$A^1E' = E'\times_{\tE_X} M^1\tE_X ,$$
  and  hence that
  $$A^1E'_X = E'_X\times_{\tE_X} M^1\tE_X.$$
  Since
  $   E'_X\times_{\tE_X} M^1\tE_X= (M^1 \tE_X)^\nabla$ and 
since the $p$-curvature of $M^1\tE_X$ vanishes, it follows
  by Cartier descent that the natural map
    $F^*(A^1E'_X) \to M^1\tE_X$ is an isomorphism.
Then statement (2) for  $i =0$ implies that the map
  $F^*\gr^0_A E' \to \gr^0_M \tE$ is an isomorphism.
  By Cartier descent again, it follows that  the map
  $\gr^0_A E' \to (\gr^0_M \tE)^\tnabla$ is an isomorphism.
  This proves (3).  Since $\gr^0_A E' = \gr^0_A E'_X$ and
  $\gr^0_M \tE = \gr^0_M \tE_X$, (4) follows as well.

  For the induction step, assume that  the statements hold for  all $j
  < i$. 
 We have a commutative diagram with exact rows:
  \begin{diagram}
    0 & \rTo&   F^* A^{i}E' & \rTo & F^* A^{i-1} E' & \rTo & F^*\gr^{i-1}_A E' & \rTo&0\cr
      && \dTo && \dTo && \dTo \cr
 0 &\rTo&   M^{i} \tE & \rTo & M^{i-1} \tE & \rTo &\gr^{i-1}_M\tE & \rTo &0\cr
  \end{diagram}
Statement (1) for $i -1$ implies that the middle vertical arrow is an
isomorphism and statement (4) for $i-1$ implies that the right vertical arrow is an
isomorphism.  Thus the left vertical arrow is also an isomorphism,
proving statement (1) for $i$.

We also  have a commutative diagram
with exact rows:
  \begin{diagram}
    0 & \rTo&   F^* A^{i}E'_X & \rTo & F^* A^{i-1} E'_X & \rTo & F^*\gr^{i-1}_A E'_X & \rTo&0\cr
      && \dTo && \dTo && \dTo \cr
 0 &\rTo&   M^{i} \tE_X & \rTo & M^{i-1} \tE_X & \rTo &\gr^{i-1}_M\tE_X & \rTo &0\cr
  \end{diagram}
  Statements (2) and (3) for $i-1$ imply that
  the  two vertical  maps  on the right are isomorphisms, and
  consequently so is the map on the left.
  Since the $p$-curvature of $M^i\tE_X$ vanishes,  statement (2) holds for $i$.  

  The main difficulty is in the following lemma, which corresponds
  to \cite{O0} 2.2.1 and will allow us to prove statement (3).
  \begin{lem}\label{str.l}
   If  the statements of  Theorem~\ref{fspan.t} hold for all $j < i$, then
the map   $\gr^i_A E'_X \to F_* \gr^i_M \tE_X$ is injective.
  \end{lem}
  \begin{pf}
    Suppose $a' \in A^iE'$ lifts an element of the kernel of the map
    in the lemma.
    Then $\eta(a') = pb  + c$, with $c \in M^{i+1}\tE$ and $b \in \tE$.
    Suppose that in fact
    $$\eta(a') = p^jb +c$$
    with $c \in M^{i+1}\tE$ and $j > 0$.
    Since $\eta(a') \in M^iE$, it follows that $b \in M^{i-j}\tE$.
    Note that if $j > i$, in fact $a' \in A^{i+1}E'$ and we are done.
    On the other hand, if $0 < j \le i$,     we calculate:
      \begin{eqnarray*}
    \tPhi(\eta(a')) & = & p^j\tPhi(b) + \tPhi(c) \\
\tPhi_{i}(\eta(a')) &= &\tPhi_{i-j}(b) +  p\tPhi_{i+1}(c)\cr
\nabla \tPhi_{i}(\eta(a')) &= &\nabla \tPhi_{i-j}(b) +  p\nabla\tPhi_{i+1}(c)\cr
      \end{eqnarray*}
      Since $\nabla \tPhi_i (\eta(a')) = \tPhi_i(\tnabla(\eta(a'))$ is
      divisible by $p$, the same is true of
$\nabla\tPhi_{i-j}(b)$.
We saw in  Proposition~\ref{agft.p} that the map $\tPhi_{i-j}$
induces a horizontal isomorphism
  $\gr^{i-j}_M \tE_X \cong \gr^{j-i}_NE_X$,  and we conclude
  that the image of $b$ in $\gr^{i-j}_M \tE_X$ is horizontal.
  Statement (3) for $i-j$ says that 
  $\gr^{i-j}_A E'_X  \cong (\gr^{i-j}_M \tE_X)^\tnabla$, so there exist $b' \in A^{i-j}E',
  b'' \in \tE$, and $b''' \in M^{i-j-1}\tE$  such that
  $$b = \eta(b') + b'' + pb'''.$$
  Then
  $$ \eta(a' -p^jb') = p^jb +c -p^j\eta(b') = p^jb'' + p^{j+1}b''' +c. $$
  Since $p^jb'' \in M^{i+1}\tE$, we have can take $a'' :=
  a'-p^jb'$ and $c' = p^jb'' +c$, so now
  $\eta(a'') = p^{j+1}b''' + c'$.  Continuing by induction,
  we see that eventually
$a'$ may be chosen to lie in $A^{i+1}E'_X$.  
  \end{pf}
  Lemma~\ref{str.l} implies that $A^{i+1} E'_X = A^i E'_X\times_{M^i\tE_X} M^{i+1}\tE_X$. 
  Since $A^iE'_X = (M^i\tE_X)^\tnabla$, it follows that
  $A^{i+1} E'_X =( M^{i+1}\tE_X)^\tnabla$, and since the $p$-curvature
  of $\tnabla$ on $M^{i+1}\tE_X$ vanishes, the map
  $F^*( A^{i+1}E'_X )\to M^{i+1}\tE_X$ is an isomorphism.
  Thus (2) holds for $i+1$ and (3) holds for $i$.

Finally, we have a commutative diagram with exact rows:
  \begin{diagram}
    0 & \rTo &  F^*\gr^{i-1}_A E' & \rTo^{[p]} &F^*\gr^i_A \tE & \rTo& F^*\gr^i_A E'_X & \rTo&0\cr
      && \dTo && \dTo && \dTo \cr
    0 & \rTo&   \gr^{i-1}_M\tE & \rTo^{[p]} &  \gr^i_M \tE & \rTo & \gr^i_M \tE_X & \rTo&0\cr
  \end{diagram}
  The left vertical arrow is an isomorphism by the induction
  assumption and the right vertical arrow is an isomorphism by (3).  
  It follows that the middle arrow is an isomorphism, proving (4) and completing
  the proof of the theorem.
\end{pf}

\begin{pf*}{Proofs of Theorem~\ref{desce.t}}
  First suppose that there is a  natural number $n$
  such that $p^n\tE \subseteq E$.  Then
multiplication by $p^n$ induces a horizontal map
$\tPhi \colon \tE \to E$, defining an $F$-span, and
$E = M^n\tE$.  Since $A^n E' = \eta^{-1}(M^n\tE)$, statement (1)
of  Theorem~\ref{fspan.t} tells us that the natural map
$F^*(\eta^{-1}(E)) \to E$ is an isomorphism.

For the general case,
 let $E_n := E + p^n\tE$ for each $n$.
Then $E_n$ is also invariant under $\tnabla$.  If
$E'_n := \eta^{-1}(E_n)$ the previous paragraph tells us that  the natural map
$F^*E'_n \to E_n$ is an isomorphism.
Let $E'' := \eta^{-1}(E)$.  Then
$E'/E'' \subseteq \tE/E$, and the Artin-Rees lemma
implies  that for some $r > 0$,
$(p^{n+r}\tE/E )\cap( E'/E'') \subseteq p^nE'/E'''$ for all  $n \ge
0$. Then $E_{n+r} \cap E' \subseteq E'' + p^nE'$ for all $n \ge 0$.  Taking the intersection
over  $n$, we find that
$$E'' \subseteq \cap \{ E'_n : n \ge 0 \} = \cap \{  E' \cap E_n:  n \ge 0 \}
\subseteq   \cap \{ E''+ p^nE' : n \ge 0\} = E'' .$$
Since $F$ is finite and flat,
the natural map
\[F^* E'' = F^* (\cap \{ E'_n  : n\ge 0\}) \to \cap \{ F^* E'_n : n \ge 0\} = \cap \{ E+ p^n\tE : n \ge 0\} = E\]
is an isomorphism.

Let us now  review Shiho's Theorem 3.1 (\cite{Sh}) and explain how it implies Theorem~\ref{desce.t}.
Shiho's theorem shows that Frobenius pullback defines an equivalence $C_F$ from the category
of  nilpotent $p$-connections on $Z'/T$ to the category of nilpotent connections on $Z/T$. 
Recall that there is a unique map
$$\zeta_F \colon  \Omega^1_{Z'/S} \to F_*\Omega^1_{Z/S}$$
such that $p\zeta_F$ is the differential of $F$.  Then if  $(E',\theta')$
is a module with $p$-connection on $Z'$, it is easy to verify that
there is a unique connection $\tilde \nabla$ on $\tilde E := F^*E'$
such that $\tilde \nabla\circ  \eta = (\zeta_F \ot \eta)  \circ \theta'$.
Shiho proves that  the functor $C_F$ taking $(E',\theta')$ to $(\tilde E, \tilde \nabla)$
is an equivalence by studying the descent data for the PD-thickenings which
correspond to the crystalline interpretations of these categories.

To apply Shiho's result, suppose that $(E',\nabla')$ is a module with nilpotent
connection on $Z'/T$.  Then $\theta':= p\nabla'$ is a nilpotent $p$-connection on
$E'$, and
\[ (\zeta_F \ot \eta)  \circ \theta' = (p\zeta_F \ot \eta)  \circ   \nabla'  = (F^* \ot \eta)  \circ \theta',\] 
and so the connection $\tnabla$ in Shiho's correspondence is just
the Frobenius pullback connection on $\tE :=F^*(E')$. Let  $E \subseteq \tE$ be a submodule
which is invariant under $\tnabla$.  We claim that the induced
connection on $E$ is also nilpotent. To see this, let $N^iE := p^i\tE \cap E$, which is also
invariant under $\tnabla$, and note that the inclusion map induces an injection
$\gr^i_NE \to p^i\tE/p^{i+1}\tE \cong \tE/p\tE$.  It follows that the $p$-curvature
of each $\gr^i_N E$ vanishes.  On the other hand,
we have a surjection $\gr^i_NE \to  (N^iE +pE)/(N^{i+1} + pE)$, so the $p$-curvature
of each of the latter also vanishes.   By Artin-Rees, there is a natural number $r$ such
that $p^r \tE \cap E \subseteq pE$. Thus the images of $N^iE$ in $E/pE$ define a finite
and exhaustive filtration of $E/pE$, and so the connection on $E$ is
indeed nilpotent. 
 Then the  full faithfulness of $C_F$
implies that there is a submodule $E''$ of $E'$, stable under $\theta' := p\nabla$,
such that $F^*E'' = E \subseteq \tE$.  Necessarily $E' \subseteq\eta^{-1}(E)$, and 
since $F^*E'  \cong E$, it follows that $E' = \eta^{-1}(E)$, and  the proof is complete. 
\end{pf*}

\begin{cor}\label{fdes.p}
  With the notations above, the natural maps
  \begin{eqnarray*}
    F^*(A^iE') & \to&  M^i  \tE \cr
    F^*(A^{[i]} E') & \to&  M^{[i]}  \tE \cr
    A^{i}  E' &\to &\eta^{-1} (M^{i}  \tE) \cr
      A^{[i]}  E' &\to &\eta^{-1} (M^{[i]}  \tE )\cr
      \end{eqnarray*}
are isomorphisms. 
\end{cor}
\begin{pf}
  The following lemma shows that the first pair of
  equations of the proposition implies  the second  pair.
  (Note:  The third equation is true by definition,
  but we shall need the fourth equation, which does
  not seem to be so obvious.)

\begin{lem}\label{back.l}
  If $A $ is an $\cO_{X'}$-submodule of $ E'$ and $M$ is the image
of $F^*A $ in $\tE'$, then $A = \eta^{-1}_F (M)$. 
\end{lem}

 \begin{pf}
    Since $F$ is flat, the map $F^*A \to F^*E'$ is injective. 
    Let $A' := \eta^{-1}_F(M)$.  Then $A'$ is an $\cO_{X'}$-submodule  of $E'$ and $A \subseteq A'$.  We have injections
    $F^*A \to F^*A' \to M$ whose composition is an isomorphism.
Then the map $F^*A \to F^*A'$
    is an isomorphism, and   since $F$ is faithfully flat, $A = A'$.    
\end{pf}

Statement (1) of Theorem~\ref{fspan.t} proves the first equation of Proposition~\ref{fdes.p}.
Since $F$ is flat, the natural maps
  $F^*A^{i}_FE' \to F^*E'$ and 
  $F^*A^{[i]}_FE' \to F^*E'$ are  injective, 
  and moreover
$$ F^*A^{[i]}_FE' = \sum_j p^{[i-j]}F^*A^{i-j}_FE' =M^{[i]}_F F^*E' $$
This is the second equation of the proposition, which implies the
fourth, by Lemma~\ref{back.l}.
\end{pf}

If $G$ is another lifting of $F_X$, the connection
$\nabla'$ furnishes a horizontal isomorphism
$$\ep_{G,F} \colon G^*(E',\nabla') \to F^*(E',\nabla').$$
and hence a map
$$\tilde \Phi_G\colon  G^*(E',\nabla') \to (E,\nabla),$$
making the diagram
\begin{diagram}
  G^*E' &\rTo^{\tilde \Phi_G} & E \cr
  \dTo^\ep & \ruTo_{\tilde \Phi_F} \cr
   F^*E'
\end{diagram}
 commutative.  It follows 
that the isomorphism $\ep$ takes
$M^i_G G^*E'$ isomorphically to $M^i_F F^*E'$,
and that $N^{-i}_F E = N^{-i}_GE$. 

The diagram:
 \begin{diagram}
    E' & \rTo^{\eta_G}  & G_* G^*E' & \rTo&G_*E \cr
                                 && \dTo^\ep  &\ruTo_{G_*\tilde \Phi_F}\cr
                                   && G_* F^*E'
  \end{diagram}
  shows that
  \begin{equation}
    \label{eq:4}
    \Phi_G =G_*(\tilde\Phi_F) \circ \ep_{G,F}.
  \end{equation}

  It need not be the case that
  $A^i_F E' = A^i_G E'$; we give an example below.
  The following proposition remedies this situation.

 \begin{prop}\label{main.p}
   Let $F$ and $G$ be liftings of $F_{X/S}$ to maps $Z \to Z'$.  
   Then for all $i$,  $A^{[i]}_F E' = A^{[i]}_G E'$.
 \end{prop}
 \begin{pf}
 
  For simplicity we write the rest of the proof assuming
  that $X$ is a curve and that
it admits a local coordinate $t$.
  Let $\partial := \nabla'(d/dt)$ and suppose that
  $G^*(t) = F^*(t) + pg$.  
  Then if $e' \in E'$,
  \begin{equation}
    \label{eq:2}
    \ep(\eta_G (e')) = \sum_j p^{[j]}g^j \eta_F(\partial^j(e'))
  \end{equation}
  Note that, since the connection $\nabla'$ is nilpotent,
  the sequence $\partial^j(e')$ converges to zero.

  To prove that each $A^{[i]}_GE' \subseteq A^{[i]}_FE'$,
  it will suffice to prove that each
  $A^{i}_GE' \subseteq A^{[i]}_FE'$,
  We work by induction on $i$.
  Assume that $e' \in A^i_GE'$.  Then
  $\eta_G(e') \in M^i_G G^*E'$ and hence
  $\ep (\eta_G(e')) \in M^i_F F^*E'$.
  Then
  $$\eta_F(e') = \ep(\eta_G(e')) - \sum_{j>0}
  p^{[j]}g^j\eta_F(\partial^j(e')).$$
  By Griffiths transversality, $\partial^j(e') \in A^{i-j}_GE'$,
 hence the  induction hypothesis implies that 
 $\partial^j(e') \in A^{[i-j]}_FE'$ when $j >0$.   Then
 $p^{[j]}\partial^j(e') \in A^{[i]}_FE'$, and we conclude
 that $\eta_F(e') \in F^* (A^{[i]}_FE') + M^i_FF^*E' = M^{[i]}_F
 F^*E'$
 by  the second equation of  Proposition~\ref{fdes.p}.
Its  last equation then implies
 that $e' \in A^{[i]}_F E'$. 
\end{pf}



  
{\bf Example:} 
  Let $X := \Spec k[t,t^{-1}]$,  let $Z$ be the formal
  completion of $\Spec W[t, t^{-1}]$ and let $F$ send
  $t$ to $t^p$.  Let $E$ be the free $\cO_Z$-module
  with basis $(e_0, \ldots, e_{p-1})$, let $\eta$
  be the endomorphsm of $E$ sending $e_i$ to $e_{i-1}$
  and $e_0$ to zero, and let
  $$\nabla e_i = \eta(e_i) dt/t.$$
  Recall that
  $$\log (1+x) = x -x^2/2 +x^3/3 + \cdots.$$
  For each $i$, we have a formal horizontal section
  \begin{eqnarray*}
\tilde e_i& :=&e^{(-\log t) \eta} (e_i)   \\
           & = & e_i - (\log t) e_{i-1} + (1/2) (\log t)^2 e_{i-2} + \cdots
  + ((-1)^{i}/i!) (\log t)^i e_0.
      \end{eqnarray*}

Then
  \begin{eqnarray*}
\ e_i& =&e^{(\log t) \eta} (\tilde e_i)   \\
           & = & \tilde e_i - (\log t)\tilde  e_{i-1} + (1/2) (\log t)^2 \tilde e_{i-2} + \cdots  + ((-1)^{i}/i!) (\log t)^i \tilde e_0.
      \end{eqnarray*}
 
  Let $(E',\nabla') = (E,\nabla)$ and define
  $$\tilde \Phi_F \colon F^*(E',\nabla') \to (E,\nabla)
  : e'_i \mapsto p^i e_i.$$
  It is immediate to check that this map is horizontal

  Now suppose that $G$ is another lift of $F_X$,
  sending $t$ to $t^p + pg$.
  Let $u := (1+ pt^{-p}g)$, so that
  $t^p + pg = t^pu$. Note that
  $\log u \in p W\{t, t^{-1}\}$,
  say $\log u= p \delta$.
  That is,
  $$\delta = t^{-p}g -(p/2)t^{-2p}g^2 + (p^2/3)t^{-3p}g^3 + \cdots.$$

  To calculate $\ep := \ep_{G,F}$, we use the fact that
  it acts as the identity on horizontal sections.
  Thus
  \begin{eqnarray*}
    \ep(G^* (e_i)) & =& \ep( G^*(e^{(\log t)\eta} \tilde e_i )\cr
                              & = &  e^{\log (t^p u)\eta} \ep (G^*(\tilde  e_i) ) \cr
                       & = &  e^{\log (t^pu)\eta} F^*(\tilde e_i) \cr
               & = & e^{\log (t^pu)\eta} F^* (e^{-(\log t)\eta} )(e_i) \cr
      & = & e^{\log (t^pu)\eta}  (e^{-(\log t^p)\eta} )(e_i) \cr
       & = & e^{(\log u)\eta} e_i \cr
     & = & e_i + (\log u) e_{i-1} + (\log u)^2/(2!) e_{i-2}+ \cdots +(\log u)^i/i! e_0 \cr
     & = & e_i +  p\delta  e_{i-1} + p^{[2]}\delta^2 e_{i-2}+ \cdots +p^{[i]}\delta^i e_0 \cr
  \end{eqnarray*}
  Then 
  \begin{eqnarray*}
     \Phi_G (e_i) &=&  \tilde \Phi_F(\ep (G^*(e_i))) \\
             & = & = p^i( e_i  +  \delta e_{i-1} + \delta^{[2]}e_{i-2}
                   + \cdots + \delta^{[i]} e_0)
  \end{eqnarray*}
  In particular,
  $e_i \in A^i_G E'$ if $i < p$, but
  $\Phi_G(e_p) = p^p( \cdots) + (p^{p-1}/(p-1)!) \delta^pe_0$.

  Take for example $p = 2$.
  Then
  \begin{eqnarray*}
    \Phi_G(e_0) & =& e_0 \\
    \Phi_G(e_1) & = & 2e_1 + 2\delta e_0\\
    \Phi_G(e_2) &= &4 e_2 + 4\delta e_1 + 2\delta^2e_0
  \end{eqnarray*}
  Assume $\delta$ is not zero modulo $2$ (e.g.
  if $g = 4t^2, u = 5$, and $\delta = 1/2 \log 5$.
  Then $e_2'$ does not belong to $A^2_GE'$,
  although it does belong to $A^2_FE'$.
  On the other hand, since $G$ lifts $F_X$ and $k$ is perfect,
  there exists an element $\delta'$ of $W\{t, t^{-1}\}$ such
  that $G^*(\delta') \equiv \delta^2\pmod 2$.
Then $e'_2 -2\delta' e'_0 \in A^2_G E'$.

\end{document}